\documentclass[12pt]{amsart}
\usepackage{newcent}
\usepackage{xparse,xfrac}
\usepackage[utf8]{inputenc}
\NewDocumentCommand{\tens}{t_}
 {
  \IfBooleanTF{#1}
   {\tensop}
   {\otimes}
 }
\NewDocumentCommand{\tensop}{m}
 {
  \mathbin{\mathop{\otimes}\displaylimits_{#1}}
 }
\usepackage{amsmath,amssymb,amsthm,soul,color,framed,mathtools}
\usepackage[hyphens]{url}
\usepackage{hyperref}
\usepackage[hyphenbreaks]{breakurl}
\usepackage[margin=1in]{geometry}
\usepackage{enumerate}
\usepackage{mathrsfs}
\usepackage{tikz-cd}
\tikzset{Isom/.style={above,every to/.append style={edge node={node [sloped, auto=false]{$\sim$}}}}}
\usepackage[backend=biber,style=alphabetic,
citestyle=alphabetic]{biblatex}
\usetikzlibrary{decorations.pathmorphing}

\usepackage[cal=euler,calscaled=1.03]{mathalpha}

\usepackage{thm-restate}

\usepackage{xspace}
\xspaceaddexceptions{]\}}

\title{Higher-weight Jacobians}

\author{Sheela Devadas}
    \address{University of Maine, Orono, ME, USA}
\email{sheela.devadas@maine.edu}
    \author{Max Lieblich}
    \address{University of Washington, Seattle, WA, USA}
    \email{lieblich@uw.edu}
\date{August 2026}
\newcommand\simto{\stackrel{\textstyle\sim}{\smash{\longrightarrow}\rule{0pt}{0.4ex}}}
\newcommand{\eps}{\epsilon}
\renewcommand{\O}{\mathscr{O}}
\newcommand{\Or}{\mathcal{O}}
\newcommand{\OA}{\Or(\A)}
\newcommand{\OK}{\Or_K}
\newcommand{\Z}{\mathbf{Z}}
\newcommand{\C}{\mathbf{C}}
\newcommand{\Q}{\mathbf{Q}}

\newcommand{\Con}{\mathfrak{C}}
\newcommand{\Lat}{\mathfrak{L}}
\newcommand{\s}{\sqrt{-3}}
\newcommand{\R}{\mathbf{R}}
\newcommand{\F}{\mathbf{F}}
\newcommand{\G}{\mathbf{G}}

\newcommand{\A}{\mc{A}}
\newcommand{\Hom}{{\rm Hom}}
\newcommand{\End}{{\rm End}}
\newcommand{\Gal}{{\rm Gal}}
\newcommand{\lcm}{{\rm lcm}}
\newcommand{\coker}{{\rm coker}}
\renewcommand{\H}{{\rm H}}

\newcommand{\im}{{\rm im}}
\newcommand{\ol}[1]{\overline{#1}}
\newcommand{\mc}[1]{\mathcal{#1}}

\newcommand{\ms}[1]{\mathscr{#1}}
\DeclareMathOperator{\Cl}{Cl}

\newcommand{\lr}[1]{\langle #1\rangle}
\newcommand{\surj}{\twoheadrightarrow}

\newcommand{\smat}[1]{\left(\begin{smallmatrix} #1  \end{smallmatrix}\right)}

\newtheorem{theorem}{Theorem}[subsection]
\renewcommand{\thetheorem}{%
  \ifnum\value{subsection}>0
    \thesubsection
  \else
    \thesection
  \fi
  .\arabic{theorem}%
}
\newtheorem{corollary}[theorem]{Corollary}
\newtheorem{proposition}[theorem]{Proposition}
\newtheorem{lemma}[theorem]{Lemma}
\theoremstyle{definition}
\newtheorem{notation}[theorem]{Notation}
\newtheorem{observation}[theorem]{Observation}
\newtheorem{question}[theorem]{Question}
\newtheorem{definition}[theorem]{Definition}
\newtheorem{remark}[theorem]{Remark}
\newtheorem{examplex}[theorem]{Example}
\newenvironment{example}
  {\pushQED{\qed}\examplex}
  {\popQED\endexamplex}

\DeclareMathOperator{\Br}{Br}

\DeclareMathOperator{\rk}{rk}

\newcommand{\Zsi}{\mathfrak{R}'}
\newcommand{\Zti}{\mathfrak{R}}
\newcommand{\Qz}{\mathfrak{K}}
\newcommand{\Zfe}{\mathfrak{O}'}
\newcommand{\Ztws}{\mathfrak{O}}
\newcommand{\la}{\mathcal{L}}

\DeclareMathOperator{\Kum}{\mathcal{K}}

\begin{document}

\begin{abstract}
    We define and study Jacobians of Hodge structures with weight greater than 1.
    A complex variety has a Jacobian of weight 2 precisely when it has maximal Picard number; these weight 2 Jacobians
  naturally arise in the context of the Brauer group and the Tate conjecture. The case of surfaces of maximal Picard number has been previously studied by Beauville; the higher-dimensional case is also related to the work of Totaro on Hodge structures with no middle pieces.

  Higher-weight Jacobians are complex tori, and it is generally quite difficult to tell if they are algebraic. However, for abelian varieties of maximal Picard number, we are able to explicitly calculate their higher-weight Jacobians using algebraic number theory, and prove that they are not just abelian varieties but in fact also have maximal Picard number. We also compute higher-weight Jacobians for Kummer varieties and singular K3 surfaces. Via class field theory, we study the field of definition of abelian surfaces and singular K3 surfaces using their weight 2 Jacobians.
  

\end{abstract}
\maketitle

\vspace{-2\baselineskip}

\tableofcontents
\section{Introduction}

Classically, one can realize the Jacobian of a complex genus $g$ curve $C$ as a complex torus using the exponential sequence:
$$\H^0(C,\ms O)\to\H^0(C,\G_m)\to\H^1(C,\Z(1))\to\H^1(C,\ms O)\to\H^1(C,\G_m)\to\H^2(C,\Z(1)).$$ We know that the complex span of the image of $\H^1(C,\Z(1))$ is all of $\H^1(C,\ms O)$; the latter is a complex vector space of dimension $g$. The image is discrete (a free abelian group of rank $2g$) and the cokernel of this map is a complex torus. 

If we carry out the same analysis for a smooth projective complex variety $X$ of arbitrary dimension, we get the Picard variety of $X$: the complex span of the integral cohomology is always the full coherent cohomology (by Hodge theory). The deepest subtlety in this classical setting is the \emph{algebraicity\/} of the resulting complex torus; since the torus admits a polarization, it is in fact not just a complex torus but an abelian variety. 


The question then arises in higher cohomological degrees: when can one attach an abelian variety to $X$ that measures the connected part of a cohomology group? To be specific, degree 2 is of interest because of the Brauer group $\H^2_{\textrm{\'et}}(X,\G_m)$. However, what draws interest to this case is also an immediate roadblock: for smooth $X$, the Brauer group $\H^2(X,\G_m)$ (in fact, all higher degree cohomology) is always torsion. Moreover, the map $\H^2(X,\Z(1))\to\H^2(X,\ms O)$ need not have discrete image; that is only true when the Picard number of $X$ is maximal. Finally, the cohomology has even degree, so there is no reason for any torus arising in this construction to be algebraic. Indeed Weil's \cite{weil} and Griffith's \cite{griff, griff2} notions of an intermediate Jacobian, although constructed differently by putting a complex structure on the real torus $\H^n(X,\R)/\H^n(X,\Z)$, are only defined in odd degree. 

We will then define a ``higher-weight Jacobian'' which generalizes the complex torus associated to $\H^2$ for a surface of maximal Picard number discussed in \cite{beauville}; the more specific weight $2$ cases of singular abelian surfaces \cite{sm74} or singular K3 surfaces \cite{si77} have also been considered previously. Higher-dimensional abelian varieties of maximal Picard number have been studied in \cite{kollar2003}; see also the related endpoint Hodge structures $(n, 0, \dots, 0, n)$ whose endomorphisms are studied in \cite{totaro}.

\begin{notation}\label{notation:latcon}
    Given a (pure) integral Hodge structure $H$ of weight $m$ such that $H^{i,j}=0$ for $i<0$ or $j<0$, define $\Lat(H)$ to be the image of the natural map $H\to H^{0,m}$ and let $\Con(H)= H^{0,m}/\Lat(H)$. If $H=\H^m(X,\Z)$ for a smooth complex projective variety $X$, we will write $\Lat_m(X)=\Lat(H)$ and $\Con_m(X)=\Con(H)$. In this case, $$\Con_m(X)=\H^m(X,\ms O_X)/\Lat_m(X)=\coker(\iota_m),$$
    where $\iota_m:\H^m(X,\Z)\to\H^m(X,\ms O_X)$ is the natural map induced by functoriality.
\end{notation}

\begin{definition} 
    We see that a pure integral Hodge structure $H$ of weight $m > 0$ \emph{has a Jacobian\/} if $\Lat(H)$ is a discrete subgroup of the complex vector space $H^{0,m}$. This gives $\Con(H)=H^{0,m}/\Lat(H)$ a natural structure of a compact complex torus $\C^g/\Gamma$.\footnote{Throughout this paper, we use ``complex torus'' to refer to a compact complex torus.} 
    If $H=\H^m(X,\Z)$ with $X$ a smooth complex projective variety, we instead say that $X$ \emph{has a Jacobian of weight $m$}, or an \emph{$m$-Jacobian}. 
    When $m=2$, we will also use the name \emph{Brauer-Jacobian\/}, and we will write $\Con$ for $\Con_2$.
\end{definition}
   \begin{remark}\label{rmk:rhomaxl} We note that for a smooth projective complex variety $X$, the cokernel of the map $\H^2(X,\Z) \to \H^2(X,\O)$ is a complex torus when the Neron-Severi group ${\rm NS}(X)$ is of ``maximal'' rank; in \cite{beauville} this is referred to as $X$ being \emph{$\rho$-maximal}. In particular a complex abelian variety is $\rho$-maximal if and only if it is isogenous to a self-product of a CM elliptic curve.

In fact for complex abelian varieties which are $\rho$-maximal, 
not only is the cokernel of the map $\H^2(X,\Z) \to \H^2(X,\O)$ a complex torus, we will show (Lemma~\ref{lem:definingBmX}) that the same is true of the map $\H^m(X,\Z) \to \H^m(X,\O)$. 
We will compute $\Con_m(X)$ for $\rho$-maximal $X$ and $2 \le m \le \dim(X)$ in Theorem~\ref{thm:computingBmX}. 
\end{remark}

\begin{definition} Similarly to the term $\rho$-maximal of \cite{beauville}, we say a complex variety $X$ is \emph{$q$-maximal} if the cokernel of the map $\H^q(X,\Z) \to \H^q(X,\O)$ is a complex torus, so $X$ has a $q$-Jacobian. We see that $X$ is $2$-maximal if and only if it is $\rho$-maximal if and only if it has a Brauer-Jacobian. 
\end{definition}

In particular, when $C$ is a curve, the ``weight $1$ Jacobian of $C$'' defined above coincides with the realization of the Jacobian of $C$ as a complex torus via the exponential sequence.  

Furthermore, note that there is a natural map
 $$\ker(\H^m_{\textrm{\'et}}(X,\G_m)\to\H^{m+1}(X,\Z(1)))\to\Con_m(X)$$
 whose image is the torsion subgroup of $\Con_m(X)$. That is, the topologically trivial degree $m$ cohomology is closely related to the torsion points of the $m$-Jacobian.

 In particular, when $X$ has a $2$-Jacobian or Brauer-Jacobian, the torsion points of the Brauer-Jacobian are identified with the topologically trivial part of the Brauer group $\H^2_{\textrm{\'et}}(X,\G_m)$ (whence the name). In this paper, we study examples of Brauer-Jacobians for 
several types of varieties, including CM abelian varieties and K3 surfaces, for which the topologically trivial part of the Brauer group is in fact the entire Brauer group.

We prove that the Brauer-Jacobian of an abelian variety $X$ is algebraic whenever it exists by using the product decomposition of $2$-maximal abelian varieties proved in \cite{schoen} to carry out explicit number theoretic calculations.

\begin{restatable*}{theorem}{main}\label{thm:main}
For a complex torus $X$ with $\dim X > 1$, the following are equivalent: 
\begin{enumerate}[(1)]
\item $X$ is isomorphic to some product $E_1 \times E_2 \times \dots \times E_n$ for $E_i$ pairwise isogenous CM elliptic curves;
\item $X$ is isogenous to $E^n$ for a CM elliptic curve $E$;
\item $\Con(X)$ exists as a complex torus;
\item $\Con(X)$ exists as a complex torus and in fact is an abelian variety;
\item $X$ has maximal Picard number.
\end{enumerate}
\end{restatable*}

Moreover, using  the classification of singular K3 surfaces $Y$ by quadratic forms \cite{si77}, we can also explicitly compute their Brauer-Jacobians via the isomorphism between the form class group and the ideal class group \cite[Theorem 7.7]{cox11}.

\newcommand{\thmkthreepream}{Let $Y$ be a singular K3 surface with associated quadratic form $Q_Y=ax^2+bxy+cy^2$; then $m=\gcd(a,b,c)$ is the degree of primitivity of $Q_Y$ and $\widetilde{Q_Y}=\frac{1}{m}Q_Y$ is a primitive quadratic form. Set $D={\rm disc}(Q_Y), D_0=\frac{D}{m^2}={\rm disc}(\widetilde{Q_Y})=D_0$, and let $\Or_0$ be the unique order in an imaginary quadratic field of discriminant $D_0$.}

\begin{restatable*}{theorem}{thmkthreebrauerjac}
\label{thm:k3brauerjac}
    \thmkthreepream 
    
    Then $Y$ has a Brauer-Jacobian and under the isomorphism $\Cl(D_0) \cong \Cl(\Or)$ of Remark~\ref{rmk:classisom} we have $[\widetilde{Q_Y}] = [\Con(Y)]$.
\end{restatable*}

In higher dimensions, for the Kummer variety associated to an abelian variety we have 

\begin{restatable*}{theorem}{kummercon} \label{thm:kummercon} If $\mathcal{A}$ is a complex abelian variety with a $2m$-Jacobian $\Con_{2m}(\A)$, then the Kummer variety $\mathcal{K}(\mathcal{A})$ associated to $\mathcal{A}$ also has a $2m$-Jacobian $\Con_{2m}(\mathcal{K}(\mathcal{A}))$ which is isogenous to $\Con_{2m}(\mathcal{A})$. 
\end{restatable*}

We also demonstrate that one can often calculate the isomorphism class of $\Con(X)$. For example, for two elliptic curves $E_1,E_2$ which both have CM by the same order $\Or$, we find that $$[\Con(E_1 \times E_2)] = [E_1][E_2] \in \Cl(\Or),$$ where $[E]$ denotes the element of $\Cl(\Or)$ corresponding to the isomorphism class of an elliptic curve $E$ with CM by $\Or$. If $E_1, E_2$ are isogenous CM elliptic curves, each with CM by an order $\Or_i$ in the same imaginary quadratic number field $K$, then we show that $\Con (E_1 \times E_2)$ is an elliptic curve with CM by the order $\Or_1+\Or_2 \subseteq K$, the smallest order containing both $\Or_i$ (Theorem~\ref{thm:BE1E2}). 

More generally, for $X \cong E_1 \times \dots \times E_n$ we have $$\Con(X) \cong \prod_{i < j} \Con(E_i \times E_j),$$ so $\Con(X)$ is a product of $\binom{n}{2}$ pairwise isogenous CM elliptic curves (Theorem~\ref{thm:bvarnell}).

In fact, for $X$ satisfying the equivalent conditions of Theorem~\ref{thm:main}, not only is $\Con(X)$ an abelian variety, $\Con_m(X)$ for all $m$ between $2$ and $\dim X$ is an abelian variety which has maximal Picard number (Theorem~\ref{thm:computingBmX}).

For such varieties $X$, we end up with a complex torus $B$ and a ``natural'' map from $\Br(X)$ to the torsion points of $B$.
If $B$ is itself algebraic and sufficiently functorial, then one could use this geometric $2$-Jacobian for $\Br(X)$ to study it using the geometric and arithmetic properties of rational points on varieties. This would have many significant consequences, including applications to the Tate conjecture.

As should not be surprising, that is, sadly, not the case. 
We can give many examples of abelian surfaces $\A$ defined over a field $k$ for which the Brauer-Jacobian is an elliptic curve $E$ that does not have a model over $k$, demonstrating the ``transcendental'' nature of this construction.

 In fact, when $\A=E_1 \times E_2$ is primitive (so the endomorphism rings of the $E_i$ and $\Con(\A)$ are all equal) we can describe exactly for which cases $\A$ can be decomposed as a product of elliptic curves over the minimal field of definition of the elliptic curve $\Con(\A)$; this is closely related to the field of definition of the Kummer surface of $\A$. 

\begin{restatable*}{proposition}{decompfieldofdef}\label{prop:decompfieldofdef} Assume $\A=\C^2/\Lambda$ is a complex torus of dimension $2$ with a Brauer-Jacobian $\Con(\A)$ such that $\A$ is primitive.

There exist elliptic curves $E_1,E_2$ with $\A \simeq E_1 \times E_2$ such that $\Q(j(E_1))=\Q(j(E_2))=\Q(j(\Con(\A)))$  if and only if the group element $[\Con(\A)] \in \Cl(\OA)$ has order $\le 2$. 
\end{restatable*} 

\begin{restatable*}{proposition}
{kumfieldofdef}\label{prop:kumfieldofdef} Assume $\A=\C^2/\Lambda$ is a complex torus of dimension $2$ with a Brauer-Jacobian $\Con(\A)$ and degree of primitivity $m$. Let $\Or=\End(\Con(\A))$ and $\OA=1+m\Or$.

If the group element $[\Con(\A)] \in \Cl(\Or)$ has order $\le 2$, then the Kummer surface $\Kum(\A)$ has a model over the field $\Q(j(\OA))$. 
\end{restatable*} 

Given that the construction of $\Con_m(X)$ appears to be transcendental in nature, it is unexpected that the results of the construction are algebraic tori (both here and in \cite{beauville}). Note that there are some cases described in Remark 7.6 of \cite{shioda72} where it is not clear if the torus is algebraic. One is naturally led to the following question:

\begin{question}
    For which values of $m$ is there smooth projective $X$ admitting an $m$-Jacobian so that $\Con_m(X)$ is \emph{not\/} algebraic?
\end{question}

Our computations of the (Brauer-)Jacobian of a complex abelian variety rely on some background from number theory and complex algebraic geometry given in Section~\ref{sec:background}. More specifically, in Section~\ref{sub:latprod} we discuss the connection between complex elliptic curves and lattices in $\C$ in order to define a ``lattice product'' which computes $\Lat_2(E_1 \times E_2)$, and in Section~\ref{sub:fieldofdef} we consider the minimal field of definition of a complex elliptic curve using class field theory, which will be used for proving Proposition~\ref{prop:decompfieldofdef}. With the necessary background material on lattices and complex elliptic curves established, we can define and compute the Brauer-Jacobians for complex abelian varieties by considering them as complex tori $\C^n/\Gamma$ in Section~\ref{sub:computingcoho}.

 We can then consider the structure of higher-weight Jacobians as abelian varieties. We use the class group of the endomorphism algebra of a CM elliptic curve $E$ to compute the $m$-Jacobians for complex abelian varieties isogenous to a power of $E$, beginning with products $E_1 \times \dots \times E_n$ for elliptic curves $E_i$ each isogenous to $E$ and $m=2$ in Section~\ref{sub:descBX}. In Section~\ref{sub:haveJacobians} we show that in fact all complex abelian varieties having a Brauer-Jacobian can be decomposed as a product in this way (Theorem~\ref{thm:main}), and furthermore consider how the minimal field of definition for different product decompositions of complex abelian surfaces relates to the $m$-Jacobians. In Section~\ref{sub:xsimeqbx} we consider when an abelian variety could be isomorphic to its $m$-Jacobian by leveraging a particular decomposition into a product of elliptic curves from \cite{bf68}.

 Finally, in Section~\ref{sec:K3} we extend our discussion of higher-weight Jacobians for complex abelian varieties to complex K3 surfaces. In Section~\ref{sec:generalities} we discuss generalities on higher-weight Jacobians for blowups and equivariant $m$-Jacobians, which will be applied to computing $m$-Jacobians for Kummer varieties in Theorem~\ref{thm:kummercon}.

 Our calculations of $j$-invariants for various CM elliptic curves/lattices can be found in Appendix~\ref{app}.\newline

\subsection*{Acknowledgments}
We would like to thank Lalita Devadas, János Kollár, and Martin Olsson for helpful conversations and correspondence. Devadas was supported by NSF postdoc grant DMS-2102960. Lieblich was supported by NSF grant DMS-1901933 and NSF FRG grant DMS-2152235.

\section{Complex lattices and complex tori}\label{sec:background}

Lattices in the complex numbers are closely related to complex elliptic curves; this is discussed in various sources including \cite{silverman09}. In particular, by \cite[Corollary VI.5.1.1]{silverman09}, every complex elliptic curve $E$ is isomorphic as a complex Lie group to $\C/\Lambda$ for some lattice $\Lambda \subset \C$, which is unique up to homothety. 

Furthermore, by \cite[Theorem VI.4.1]{silverman09}, every isogeny of complex elliptic curves $\phi: E_1 \to E_2$ (with $E_i=\C/\Lambda_i$) arises from a homothety of $\Lambda_1$ with a sublattice of $\Lambda_2$---i.e., the isogeny corresponds to a nonzero complex number $\alpha$ with $\alpha\Lambda_1 \subseteq \Lambda_2$. We will refer to lattices with this property as \emph{isogenous}, as we do for the corresponding elliptic curves. Isogeny of lattices is an equivalence relation.

\begin{definition}
    A lattice $\Lambda \subseteq \C$ has \emph{complex multiplication} or \emph{CM} if its endomorphism ring $\End(\Lambda) := \left\{\alpha \in \C : \alpha\Lambda \subseteq \Lambda\right\}$ satisfies $\Z \subsetneq \End(\Lambda)$. 
\end{definition}

A discussion of orders in imaginary quadratic number fields and their relation with lattices with complex multiplication can be found in \cite[Section 7]{cox11}. It is well-known that the endomorphism ring of a complex lattice with CM is an order $\Or$ in an imaginary quadratic number field $K$ \cite[Lemma 7.2]{cox11}. Such an order is equal to $\Z+f \OK$ for $f$ the index of $\Or$ in $\OK$; this $f$ is also called the \emph{conductor} of $\Or$. It is also well-known that a lattice $\Lambda=\lr{1,\tau}$ has complex multiplication if and only if $\tau$ is the root of a polynomial $ax^2+bx+c \in \Z[x]$.

In fact a CM lattice $\Lambda$ is always homothetic to an ideal of its endomorphism ring $\Or=\End(\Lambda)$; thus we consider specifically fractional ideals of $\Or$, that is, finitely generated $\Or$-submodules of $K$, which are homothetic to an ideal of $\Or$ via ``clearing denominators''. We say that a fractional ideal $\Lambda$ of $\Or$ is \emph{invertible} if there exists another fractional ideal $\Lambda^{-1}$ such that the ideal product $\Lambda\Lambda^{-1}=\Or$. This condition is equivalent to $\End(\Lambda)=\Or$.

The set of homothety classes of lattices with endomorphism ring $\Or$ is in bijection with the ideal class group $\Cl(\Or)=\mc{I}(\Or)/\mc{P}(\Or)$ (for $\mc{I}(\Or)$ the group of invertible fractional $\Or$-ideals under multiplication, and $\mc{P}(\Or)$ the subgroup of principal fractional $\Or$-ideals), since a lattice $\Lambda$ with $\End(\Lambda)=\Or$ is homothetic to an ideal of $\Or$, and a homothety of (fractional) ideals of $\Or$ can be considered as multiplication by a principal fractional ideal.



\begin{lemma}\label{lem:isog} If $\Or_1,\Or_2$ are orders in an imaginary quadratic number field $K$ and $\Lambda_1,\Lambda_2$ are lattices with CM by $\Or_1,\Or_2$, respectively, then $\Lambda_1$ is isogenous to $\Lambda_2$. \end{lemma}
\begin{proof}
The lattice $\Lambda_1$ with CM by $\Or_1$ is homothetic to an ideal of $\Or_1$, which is a sublattice of $\Or_1$. Hence $\Lambda_1$ is isogenous to $\Or_1$ (and similarly for $\Lambda_2,\Or_2$). Furthermore, if $\Or_i$ has conductor $f_i$, then both $\Or_1,\Or_2$ are sublattices of the order in $K$ with conductor $\gcd(f_1,f_2)$, so they are also isogenous to one another. It follows that $\Lambda_1$ is isogenous to $\Lambda_2$. 
\end{proof}

\subsection{Lattice product}\label{sub:latprod}

The ``lattice product'' defined below is directly related to the Brauer-Jacobian of the product of two elliptic curves; in Lemma~\ref{lem:prod2ell} we will compute that for two isogenous CM elliptic curves $\C/\Lambda_1, \C/\Lambda_2$, the Brauer-Jacobian is isomorphic to the elliptic curve corresponding to this lattice product $\Lambda_1\Lambda_2$, which we distinguish from the direct product lattice $\Lambda_1 \times \Lambda_2 \subseteq \C^2$.

\begin{definition}\label{def:latprod}
Given two lattices $\Lambda_1,\Lambda_2 \subseteq \C$, the \emph{lattice product} $\Lambda_1  \Lambda_2$ or $\Lambda_1 \cdot \Lambda_2$ is the additive subgroup of $\C$ consisting of finite sums of products of an element of $\Lambda_1$ with an element of $\Lambda_2$, or $$\Lambda_1  \Lambda_2 = \left\{\sum_{i=1}^s a_ib_i : s \in \mathbf{N}, a_i \in \Lambda_1,b_i \in \Lambda_2\right\}.$$
\end{definition}

\begin{lemma}\label{lem:latprodproperties}
\begin{enumerate}[(i)]
\item The lattice product is commutative and associative.
\item If $\Lambda_1=\lr{v_1,w_1},\Lambda_2=\lr{v_2,w_2}$, then $\Lambda_1  \Lambda_2 = \lr{v_1v_2,v_1w_2,w_1v_2,w_1w_2}$.
\item If $\Lambda_1,\Lambda_3$ are homothetic, and $\Lambda_1 \Lambda_2$ is a lattice, then $\Lambda_1  \Lambda_2,\Lambda_2  \Lambda_3$ are homothetic and the latter is also a lattice.
\item If $\Or=\End(\Lambda)$, then $\Or\Lambda  = \Lambda$. 
\item For two orders $\Or_1,\Or_2$ of an imaginary quadratic number field $K$ with conductors $f_1,f_2$ respectively, then $\Or_1  \Or_2=\Or_3$ for $\Or_3$ the order of $K$ with conductor $f_3=\gcd(f_1,f_2)$. ($\Or_3$ can alternatively be described as $\Or_1+\Or_2$, the ``smallest'' order of $K$ that contains both $\Or_1$ and $\Or_2$.)
\end{enumerate}
\end{lemma}
\begin{proof}
Properties (i)-(iii) are all immediate consequences of the definition of the lattice product; properties (iv) and (v) are both simple computations.
\end{proof}
\begin{lemma}\label{lem:latprodislat}
If $\Lambda_1=\lr{v_1,w_1}$ and $\Lambda_2=\lr{v_2,w_2}$ are isogenous lattices with complex multiplication, then $\Lambda_1\Lambda_2=\lr{v_1v_2,v_1w_2,w_1v_2,w_1w_2}$ is also a lattice.
\end{lemma}
\begin{proof} 
After applying homotheties, we can assume that $\Lambda_1=\lr{1,\tau_1},\Lambda_2=\lr{1,\tau_2}$ for $\tau_1,\tau_2 \in \C$.

Since the $\Lambda_i$ have CM, both fields $\Q(\tau_i)$ are quadratic extensions of $\Q$. By hypothesis, $\Lambda_2$ is homothetic to a sublattice of $\Lambda_1$. 
It follows that $\Q(\tau_1) = \Q(\tau_2)$. Hence the $\Q$-vector space $\Q \otimes_\Z \Lambda_1\Lambda_2=\Q\langle 1, \tau_1, \tau_2, \tau_1\tau_2\rangle $ is 2-dimensional. As $\Lambda_1\Lambda_2 \subseteq \C$ is necessarily non-torsion, it must be a free $\Z$-module of rank $2$ and hence a lattice. \end{proof}




For ideals and fractional ideals there is also a notion of a product, which would be written the same way as the lattice product; however, there is no ambiguity of notation as we show in the following lemma.

\begin{lemma}\label{lem:latprodidealprod}
If $\Lambda_1,\Lambda_2$ are fractional ideals of the same order $\Or$, then the lattice product $\Lambda_1  \Lambda_2$ is equal to their product as fractional ideals.
\end{lemma}
\begin{proof} 
Both $\Lambda_1,\Lambda_2$ are isogenous to $\Or$ as they are homothetic to an ideal (which is a sublattice) of $\Or$. Thus they are isogenous CM lattices and their lattice product is a lattice by Lemma~\ref{lem:latprodislat}. Furthermore, the definition of a product of fractional ideals coincides with the definition of the lattice product (finite sums of products $ab$ with $a \in \Lambda_1, b \in \Lambda_2$). \end{proof}


\begin{lemma}\label{lem:latprodextscalars}
Consider two orders $\Or \subseteq \Or'$ of an imaginary quadratic field $K$, and let $\Lambda$ be an invertible fractional ideal of $\Or$. Then $\Or'\Lambda$ is an invertible fractional ideal of $\Or'$, and it is isomorphic as an $\Or'$-module to the extension of scalars $\Lambda \tens_{\Or} \Or'$.
\end{lemma}
\begin{proof}
By Lemma \ref{lem:isog}, $\Lambda$ is isogenous to $\Or$, so we know that $\Or'\Lambda$ is a lattice by Lemma \ref{lem:latprodislat}. It is also evident that $\Or'\Lambda$ is closed under multiplication by $\Or'$ and that it is contained in $K$, so $\Or'\Lambda$ is a fractional ideal of $\Or'$; we will now prove that it is invertible.

As $\Lambda$ is an invertible fractional ideal of $\Or$, it has an inverse $\Lambda^{-1}$---that is, $\Lambda^{-1}$ is another invertible fractional ideal of $\Or$ such that $\Lambda\Lambda^{-1}=\Or$, where the product can be taken either as fractional ideals or as lattices by Lemma~\ref{lem:latprodidealprod}. (In fact $\Lambda^{-1}$ is homothetic to $\ol{\Lambda}$.)

Then by Lemma~\ref{lem:latprodproperties}(i), we see that $(\Or'\Lambda)  (\Or'\Lambda^{-1})= (\Or'\Or')(\Lambda\Lambda^{-1})=\Or'  \Or  = \Or',$ and $ \Or'\Lambda^{-1} $ is a fractional ideal of $\Or'$ for the same reason that $\Or'\Lambda$ is. Hence by Lemma~\ref{lem:latprodidealprod}, we see that $\Or'\Lambda$ is an invertible fractional ideal of $\Or'$ as we have just exhibited an inverse.

There is an obvious map $\rho: \Lambda \tens_{\Or} \Or' \to \Or'\Lambda$, given by $x \tens \alpha \mapsto \alpha x$ for $x \in \Lambda,\alpha\in \Or'$. By the universal property of tensor products and the definition of $\Or'\Lambda$, this map is surjective and $\Or'$-linear. 


It remains to prove injectivity. We will use the fact that an invertible fractional ideal of an order $\Or$ is a rank $1$ projective $\Or$-module. It then follows from the above that $\Or'\Lambda$ is a rank $1$ projective $\Or'$-module, as is $\Lambda\otimes_{\Or}\Or'.$ 

For each $\mathfrak p\in \operatorname{Spec}(\Or')$, we obtain a surjective map
$(\Lambda\otimes_{\Or}\Or')_{\mathfrak p} \to 
(\Or'\Lambda)_{\mathfrak p}$. 
Both source and target are free $\Or'_{\mathfrak p}$-modules of rank $1$.
Thus, after choosing bases, we may consider $\rho_\mathfrak{p}$ as a map 
$\Or'_{\mathfrak p}\longrightarrow \Or'_{\mathfrak p}.$
Such a map is multiplication by some element of $\Or'_{\mathfrak p}$, and
surjectivity implies that this element is a unit. Therefore
$\rho_{\mathfrak p}$ is an isomorphism for every
$\mathfrak p\in \operatorname{Spec}(\Or')$.

We have now checked locally that $\rho$ is an isomorphism, so in fact 
$
\Lambda\otimes_{\Or}\Or' \cong \Or'\Lambda
$
as $\Or'$-modules.\end{proof}

\begin{lemma}\label{lem:generallatprod}
Let $\Lambda_1,\Lambda_2$ be isogenous lattices with CM, so $\Lambda_1  \Lambda_2$ is a lattice by Lemma \ref{lem:latprodislat}. Then
\begin{enumerate}[(i)]
    \item $\Or_1=\End(\Lambda_1)$ and $\Or_2=\End(\Lambda_2)$ are orders in the same imaginary quadratic number field $K$;
    \item $\End(\Lambda_1  \Lambda_2)=\Or_1  \Or_2 = \Or_3$, another order in $K$.
\end{enumerate}
\end{lemma}
\begin{proof} A ``converse'' to (i) can be found in Lemma~\ref{lem:isog}.

Up to homothety, we can assume that for some complex numbers $\tau_1,\tau_2$ we have $\Lambda_i=\lr{1,\tau_i}$; since $1 \in \Lambda_i$, we can thus assume that $\End(\Lambda_i) \subseteq \Lambda_i$ and that $\Lambda_i$ is an (invertible) fractional ideal of $\End(\Lambda_i)$. This is, in fact, a containment of lattices. Furthermore, in the proof of Lemma~\ref{lem:latprodislat}, we showed that in this situation $\tau_2 \in \Q(\tau_1)$, which is an imaginary quadratic number field because $\Lambda_1$ has CM. Then, in fact, $\Q(\tau_2)=\Q(\tau_1)=K$ is  necessarily the fraction field of both $\End(\Lambda_i)$. This proves (i).

To prove (ii), we use Lemma~\ref{lem:latprodproperties}(i) and (iv): $$\Lambda_1  \Lambda_2 = \Lambda_1  \Or_1 \cdot \Or_2  \Lambda_2 = \Lambda_1  \cdot \Or_3 \cdot \Lambda_2 = (\Or_3\Lambda_1  )  (\Or_3\Lambda_2).$$

By Lemma~\ref{lem:latprodextscalars}, both $\Or_3\Lambda_i  $ for $i=1,2$ are invertible fractional ideals of $\Or_3$; applying Lemma~\ref{lem:latprodidealprod}, we see that their lattice product is the same as their product as fractional ideals. Therefore $(\Or_3\Lambda_1  )  (\Or_3\Lambda_2)$ is also an invertible  fractional ideal of $\Or_3$, meaning $\End(\Lambda_1  \Lambda_2)=\Or_3$, as we desired to show.
\end{proof}

\subsection{Field of definition}\label{sub:fieldofdef} 

There exists a correspondence between complex elliptic curves and complex lattices, as discussed in previous sections. We note that two elliptic curves defined over number fields are isomorphic over $\C$ if and only if they have the same $j$-invariant; we generally treat the $j$-invariant of an elliptic curve and the $j$-invariant of the associated lattice as interchangeable \cite[Proposition 4.4]{silvermanAd}.

\begin{definition}\label{def:fieldofdef} If a complex elliptic curve $E$ is isomorphic to one defined over a number field $L$, then necessarily the $j$-invariant $j(E) \in L$. Conversely, as long as $j \ne 0,1728$, then $E$ is isomorphic (over $\C$) to the curve defined by $$y^2=x^3+\frac{3j}{1728-j}x +\frac{2j}{1728-j},$$ and it is easy to find examples of elliptic curves over $\Q$ with $j$-invariant $0$ or $1728$.

Therefore, the ``smallest'' number field over which an elliptic curve isomorphic to $E$ over $\C$ may be defined is $\Q(j(E))$, which we refer to as the \emph{field of definition} of $E$. \end{definition}

We will use the results of this section to relate the field of definition of a complex abelian surface $\A$ to that of its Brauer-Jacobian $\Con(\A)$, which is an elliptic curve (Proposition~\ref{prop:decompfieldofdef}).

The first main theorem of complex multiplication \cite[Theorem 11.1]{cox11} states that for an order $\Or$ of an imaginary quadratic field $K$, the $j$-invariant of any invertible fractional $\Or$-ideal $I$ (or equivalently, the $j$-invariant of any lattice with CM by $\Or$) is an algebraic integer $j(I)$, and for all invertible ideals $I \subseteq \Or$ the field extension $K(j(I))$ is the \emph{ring class field} of the order $\Or$.

The ring class field $L_\Or$ has the property that $
\Cl(\Or) \simeq \Gal(L_\Or/K)$. Under this group isomorphism, the class of an invertible fractional ideal $I$ corresponds to the element $\sigma_I$ of the Galois group which permutes $j$-invariants by $[\ol{I}]=[I^{-1}]$---i.e., $\sigma_I(j(J))=j(I^{-1}J)=j(\ol{I}J))$ for $I,J$ invertible fractional ideals of $\Or$ \cite[Corollary 11.37]{cox11}.

Furthermore, \cite[Lemma 9.3]{cox11} states that $L_\Or$ is Galois over $\Q$ with Galois group the semidirect product $\Cl(\Or) \rtimes \Z/2\Z$ where $\Z/2\Z$ acts on $\Cl(\Or)\simeq \Gal(L_{\Or}/K)$ via $\sigma \mapsto \sigma^{-1}$ and $\Z/2\Z$ acts on $L_\Or$ by complex conjugation.

This description of $\Gal(L_\Or/\Q)$ will allow us to identify when for two elliptic curves $E_1,E_2$ with CM by $\Or$ the fields of definition $\Q(j(E_1))$ and $\Q(j(E_2))$ are equal.

\begin{lemma}\label{lem:qje}
Let $E_1,E_2$ be complex elliptic curves with CM by $\Or$, an order in an imaginary quadratic field $K$. Let $L_\Or$ be the ring class field of $\Or$, which is a degree $2n$ extension of $\Q$ for $n=|\Cl(\Or)|$.
Then one of the following is true:
\begin{enumerate}[(i)]
\item $\Q(j(E_1),j(E_2))=L_\Or$, or
\item $\Q(j(E_1))=\Q(j(E_2))$ is a degree $n$ extension of $\Q$, and in $\Cl(\Or)$ the element $[E_1][E_2]^{-1}$ has order $\le 2$.
\end{enumerate}\end{lemma}
\begin{proof} If $E_1 \cong E_2$, then $j(E_1)=j(E_2)$ so case (ii) holds trivially. We may then assume below that $E_1 \not\cong E_2$. 

By \cite[Proposition 13.2]{cox11}, the $j(E_i)$ have minimal polynomial equal to the \emph{class equation} $H_\Or(X)$, which is a degree $n$ polynomial in $\Z[X]$. Hence $n=[\Q(j(E_i)):\Q]$.

It then follows from the fundamental theorem of Galois theory that the subfields $\Q(j(E_i)) \subseteq L_\Or$ are each the fixed field of an order $\frac{2n}{n}=2$ subgroup of $\Gal(L_\Or/\Q) = \Cl(\Or) \rtimes \Z/2\Z$---this order $2$ subgroup is $\Gal(L_\Or/\Q(j(E_i))) \subseteq \Gal(L_\Or/\Q)$.

Writing $\Z/2\Z$ as $\{\pm 1\}$, the order $2$ elements of $\Cl(\Or) \rtimes \Z/2\Z$ are either $(x,1)$ for $x \in \Cl(\Or)$ of order $2$, or $([E],-1)$ for any $[E] \in \Cl(\Or)$. Since $\Cl(\Or)$ acts transitively and freely on the set of $j$-invariants of elliptic curves with CM by $\Or$ \cite[Corollary 11.37]{cox11}, the subfield $\Q(j(E_1))$ cannot arise as the fixed field of $(x,1)$ for any $x \in \Cl(\Or)$ of order $2$.

Therefore for any elliptic curve $E$ with $\End(E)=\Or$, the subextension $\Q(j(E)) \subseteq L_\Or$ is the fixed field of $([I],-1)$ for some invertible fractional ideal $I$ of $\Or$. Using the description of the action of $\Cl(\Or)$ on the $j$-invariants of elliptic curves with CM by $\Or$ given in \cite[Corollary 11.37]{cox11}, we see that for any invertible fractional ideal $J$ of $\Or$, the element $([I],-1) \in \Gal(L_\Or/\Q)$ sends $j(J)$ to $j(\ol{IJ})=j(I^{-1}J^{-1})$. Therefore $([I],-1)$ fixes $j(J)$ if and only if $[J]^2=[I]^{-1} \in \Cl(\Or)$.

Hence $\Q(j(E))$ is the fixed field of $([E]^{-2},-1) \in \Gal(L_\Or/\Q)$. 

We now consider $\Q(j(E_1),j(E_2))$, which must lie between the extensions $\Q(j(E_1))$ and $L_\Or$ of $\Q$. This means $[\Q(j(E_1),j(E_2)):\Q]=n$ or $2n.$ 

In the latter case, (i) obviously holds; in the former case, $\Q(j(E_1),j(E_2))=\Q(j(E_1))=\Q(j(E_2))$, so $([E_1]^{-2},-1)=([E_2]^{-2},-1) \in \Gal(L_\Or/\Q)$. Then $[E_1]^2 = [E_2]^2 \in \Cl(\Or)$, meaning $[E_1][E_2]^{-1}$ has order $2$ (since $E_1 \not \simeq E_2$) and we have (ii) in this case, completing the proof.
\end{proof}

\begin{example}\label{ex:qjezsqrt-27} In we c
The order $\Ztws=\Z[\sqrt{-27}]$ in the field $\Qz=\Q(\sqrt{-3})=\Q(\zeta_3)$ (where $\zeta_3=\frac{-1+\sqrt{-3}}{2}$ is a cube root of unity) has ring class field $L_{\Ztws}=\Qz(\sqrt[3]{2})=\Q(\zeta_3,\sqrt[3]{2})$ and class group $\Cl(\Ztws)=\Z/3\Z$; representatives for the three homothety classes of lattices with CM by $\Ztws$ are $$\Ztws=\lr{1,3\sqrt{-3}},\qquad\la_1=\lr{3,1+\sqrt{-3}},\qquad\la_2=\lr{3,2+\sqrt{-3}},$$ the $j$-invariants of which generate the distinct degree $3$ extensions $$
    \Q(j(\Ztws))=\Q(\sqrt[3]{2}),\quad
    \Q(j(\la_1))=\Q(\zeta_3\sqrt[3]{2}),\quad
    \Q(j(\la_2))=\Q(\zeta_3^2\sqrt[3]{2}),
$$ as calculated in Example~\ref{ex:zsqrt-27} . This is in agreement with Lemma~\ref{lem:qje} as $\Cl(\Ztws)$ does not have any elements of order $2$. 
\end{example}

\begin{example}\label{ex:qjezsqrt-48}
The order $\Zfe=\Z[\sqrt{-48}]$ in the field $\Qz=\Q(\sqrt{-3})$ has ring class field $L_{\Zfe}=\Qz(i,\sqrt{2})=\Q(i,\sqrt{2},\sqrt{3})$ and class group $\Cl(\Zfe)=\Z/2\Z \times \Z/2\Z$. The four homothety classes of lattices with CM by $\Zfe$ can be computed to be represented by $$\Zfe=\lr{1,4\sqrt{-3}},\lr{4,\sqrt{-3}},\lr{4,2+\sqrt{-3}},\lr{2,1+2\sqrt{-3}},$$ the $j$-invariants of which each generate the degree $4$ extension $\Q(\sqrt{2},\sqrt{3})/\Q$ (Example~\ref{ex:zsqrt-48}). This is in agreement with Lemma~\ref{lem:qje} as every element of $\Cl(\Zfe)$ has order $2$.
\end{example}

When we consider abelian varieties that arise as the products of elliptic curves with different endomorphism rings, the following map between class groups will be useful to describe the behavior of the Brauer-Jacobian.

\begin{proposition}\label{prop:defphi}
Let $K=\Q(\alpha)$ be an imaginary quadratic field with $\OK=\Z[\alpha]$. For integers $c \mid a$, the map $\phi_{c,a,\alpha}: \Cl(\Z[a\alpha]) \to \Cl(\Z[c\alpha])$ defined by $\Lambda \mapsto \Z[c\alpha]\cdot \Lambda$ is a group homomorphism. Furthermore, for $d \mid c \mid a$ we have $\phi_{d,c,\alpha} \circ \phi_{c,a,\alpha} = \phi_{d,a,\alpha}$, and $\phi_{a,a,\alpha}$ is the identity homomorphism.
\end{proposition}
\begin{proof}
It is clear that $\phi_{c,a,\alpha}$ is well-defined on classes of invertible fractional ideals by\linebreak Lemma~\ref{lem:latprodextscalars}; in fact, it also follows that as modules over $\Z[c\alpha]$ there is an isomorphism $\Lambda \otimes_{\Z[a\alpha]}\Z[c\alpha] \cong \phi_{c,a,\alpha}(\Lambda)$. The rest of the proposition follows from Lemma~\ref{lem:latprodidealprod} and the properties of the lattice product.
\end{proof}

In the following lemma we describe the behavior of the $\phi$ map defined in Proposition~\ref{prop:defphi} with regard to $j$-invariants of class group elements.

\begin{lemma}\label{lem:qjephi}
Let $K$ be an imaginary quadratic number field with $\OK=\Z[\alpha]$, and let $c \mid a$ be positive integers. Then if $E$ is a complex elliptic curve with CM by $\Z[a\alpha]$, any elliptic curve $E'$ corresponding to the class $\phi_{c,a,\alpha}([E]) \in \Cl(\Z[c\alpha])$ has $j(E') \in \Q(j(E))$.
\end{lemma}
\newcommand{\art}{\operatorname{Art}}
\begin{proof}
Let $L_c=L_{\Z[c\alpha]},\Or_c=\Z[c\alpha]$ and define $L_a,\Or_a$ similarly. We note that the conductor of $\Or_c$ is $c$ and similarly for $\Or_a$ the conductor is $a$. 

Since $c \mid a$, we have an inclusion of ring class fields $L_c \subseteq L_a$  \cite[Exercise 9.19]{cox11} and hence a restriction map of Galois groups over $K$. 

For any order $\Or$ in $K$, the Artin symbol gives a map $\art_\Or: \Cl(\Or) \to \Gal(L_\Or/K)$. Specifically, $\art_\Or([\Lambda])=\left(\frac{L_\Or/K}{\mathfrak{p}}\right)$ where $\mathfrak{p}$ is a prime ideal of the ring of integers $\OK$ such that $\ol{\mathfrak{p}\cap \Or}=\Lambda$ and $\mathfrak{p}$ does \emph{not} divide the conductor of $\Or$ \cite[Lemma 5.19, Theorem 11.36]{cox11}. We then obtain a diagram 

\vspace{-1\baselineskip}
\[
\begin{tikzcd}\label{diag:artin}\tag{$\star$}
\Cl(\Or_a) \arrow[r, "\art_{\Or_a}"] \arrow[d, "\phi_{c,a,\alpha}"']
& \Gal(L_a/K) \arrow[d, "\mathrm{res}"] \\
\Cl(\Or_c) \arrow[r, "\art_{\Or_c}"]
& \Gal(L_c/K)
\end{tikzcd}
\]

In fact this diagram commutes. If $\Lambda$ is an invertible fractional ideal of $\Or_a$, then we may fix a prime $\mathfrak{p}$ of $\OK$ not dividing $a$ (hence not dividing $c$) with $\ol{\mathfrak{p} \cap \Or_a} =\Lambda$. From the properties of the lattice product it follows that $\ol{\mathfrak{p} \cap \Or_c}=\Or_c \cdot \Lambda = \phi_{c,a,\alpha}([\Lambda])$. Then $\mathfrak{p}$ is unramified in both $L_c$ and $L_a$ since it does not divide the conductors $c \mid a$; we may then use the definition of the Artin symbol \cite[Lemma 5.19]{cox11} (and the fact that these Galois groups and class groups are abelian) to see that $\left(\frac{L_c/K}{\mathfrak{p}}\right)=\mathrm{res}\left(\frac{L_a/K}{\mathfrak{p}}\right),$ so the diagram (\ref{diag:artin}) commutes.

 By the proof of Lemma~\ref{lem:qje} we see that $\Q(j(E))$ is the fixed field of the order $2$ subgroup of $\Gal(L_a/\Q)$ generated by $([E]^{-2},-1)$. Using the commutativity of (\ref{diag:artin}), the map $\phi_{c,a,\alpha}$ is compatible with the restriction map $\Gal(L_a/\Q) \to \Gal(L_c/\Q)$, and the image of $([E]^{-2},-1)$ under this map is $([E']^{-2},-1) \in \Cl(\Or_c) \rtimes  \Z/2\Z =\Gal(L_c/\Q)$, which we know fixes $j(E')$. Hence $([E]^{-2},-1)$ fixes $j(E')$ as well, so $j(E') \in \Q(j(E))$ as we desired to show. \end{proof}

\begin{example}\label{ex:zsqrt-36zsqrt-9}
The order $\Zti=\Z[\sqrt{-9}]=\Z[3i]$ in $K=\Q(i)$ has ring class field $L_\Zti=K(\sqrt{3})$ and class group $\Cl(\Zti)=\Z/2\Z$, and the two homothety classes of lattices with CM by $\Zti$ are represented by $\Zti=\lr{1,3i},\lr{3,1+i}$ (Example~\ref{ex:zsqrt-9}). For both these lattices, the $j$-invariant generates the quadratic extension $\Q(\sqrt{3})/\Q$.

The order $\Zsi=\Z[\sqrt{-36}]$ in the field $K=\Q(i)$ has ring class field $L_\Zsi=K(\sqrt[4]{12})$ and class group $\Cl(\Zsi)=\Z/4\Z$. The four homothety classes of lattices with CM by $\Zsi$ can be computed to be represented by $\Zsi=\lr{1,6i},\lr{3,2i},\lr{3,1\pm2i}$ (Example~\ref{ex:zsqrt-36}). The lattices $\lr{1,6i}$ and $\lr{3,2i}$ both have real $j$-invariants generating the extension $\Q(\sqrt[4]{12})/\Q$. The lattices $\lr{3, 1 \pm 2i}$ both have $j$-invariants generating the extension $\Q(i\sqrt[4]{12})/\Q$.

The map $\phi_{3,6,i}: \Cl(\Z[6i]) \to \Cl(\Z[3i])$ is the surjective map $\Z/4\Z \to \Z/2\Z$; it is easy to check that Lemma~\ref{lem:qjephi} in this case is satisfied, as $\Q(\sqrt{3}) \subseteq \Q(i\sqrt[4]{12})$ and $\Q(\sqrt{3}) \subseteq \Q(\sqrt[4]{12})$. 
\end{example}

\subsection{Computing cohomology of complex tori}\label{sub:computingcoho}

Now that we have the necessary background material on lattices and complex elliptic curves, we can define and compute the Brauer-Jacobians for complex abelian varieties that are a product of isogenous CM elliptic curves. We will compute the cohomology groups of these abelian varieties as complex tori $\C^n/\Gamma$.

The following discussion is based on Chapter I of \cite{mumabvar}.

Consider a complex torus $X=\C^n/\Gamma$, where $\Gamma$ is a lattice in $\C^n$. Let $\Gamma^*=\Hom(\Gamma,\Z)$ be the dual lattice; then we have a canonical isomorphism $\H^r(X,\Z) \simto \bigwedge^r(\Gamma^*)$ for each $r$ arising from the cup product and the natural isomorphism $\H^1(X,\Z) \simto \Gamma^*$. 

If we let $T=\Hom_\C(\C^n,\C)$ and $\ol{T}=\Hom_{\C-{\rm antilinear}}(\C^n,\C)$, then we have natural isomorphisms $\H^r(X,\O_X) \simto \bigwedge^r \ol{T}$ for each $r$, again arising from the cup product.

The natural inclusion of the constant sheaf $\underline{\Z}$ in the structure sheaf $\O_X$ gives rise to a map of the cohomology groups $\H^r(X,\Z) \to \H^r(X,\O_X)$; this map is the exterior power $\bigwedge^r$ of the map $$\Gamma^* = \Hom(\Gamma,\Z) \to \Hom(\Gamma,\Z) \otimes_\Z \C \simto \Hom_\R(\C^n,\C) \simto T \oplus \ol{T} \surj \ol{T}.$$ We will be mostly concerned with the image and cokernel of the map in the case $r=2$.

We now describe the map $\Gamma^* \to \ol{T}$ in the one-dimensional case where $\Gamma=\Lambda \subseteq \C$ is a lattice in $\C$ and $X \cong \C/\Lambda$ is a complex elliptic curve.  

We can assume that the lattice $\Lambda$ is $\lr{1,\tau}$ for some $\tau$ in the upper half plane (via homothety of the lattice). In this case, $T \simeq \C$, with the map defined by sending $\alpha \in T=\Hom_\C(\C,\C)$ to $\alpha(1)$; we see that $\alpha(x)=x\cdot\alpha(1)$ for all $x \in \C$ by the $\C$-linearity of $\alpha$. Similarly, any $\beta \in \ol{T}=\Hom_{\C-\rm{antilinear}}(\C,\C)$ satisfies $\beta(x)=\ol{x}\cdot \beta(1)$ for all $x \in \C$, so we also have an isomorphism $\ol{T} \simeq \C$ via $\beta \mapsto \beta(1)$. 

Letting $e_1,e_2$ be the $\Z$-basis of $\Lambda^*$ which is dual to the $\Z$-basis $\{1,\tau\}$ of $\Lambda$, it is a simple matrix computation to find that \begin{equation}\label{eqn:computation}\tag{$\dagger$}
    e_1 \mapsto \frac{-\tau}{\ol{\tau}-\tau},\;e_2 \mapsto \frac{1}{\ol{\tau}-\tau}
\end{equation}

via the map $\Lambda^*\to T \oplus \ol{T} \surj \ol{T} \simeq \C$ that we have described.

For a lattice $\Gamma \subseteq \C^n$ defining a complex torus $X=\C^n/\Gamma$, we see by the above that the cohomology group $\H^2(X,\Z)$ is a free abelian group of rank $\binom{2n}{2}$, and $\H^2(X,\O)$ is a complex vector space of dimension $\binom{n}{2}$. It is not hard to see that $X$ admits a Brauer-Jacobian precisely when the image of the map $\H^2(X,\Z(1))\to\H^2(X,\ms O)$ has rank $n^2-n$. 

\begin{notation}
    In the above setting, we will write $\Lat(\Gamma)$ for $\Lat(X)$, and similarly for $\Lat_m$.
\end{notation}
This permits us to think about the space of tori $X$ as parametrized by lattices $\Gamma$ in a fixed copy of $\C^n$.


\begin{lemma}\label{lem:prod2ell} If $X=E_1 \times E_2$ for isogenous CM elliptic curves $E_i=\C/\Lambda_i$, the lattice $\Lat(X)=\Lat(\Lambda_1 \times \Lambda_2)$ is homothetic to the lattice product $\Lambda_1\Lambda_2$.
\end{lemma}
\begin{proof}
    Up to isomorphism, we can assume that each $\Lambda_i=\lr{1,\tau_i}$ for $\tau_i$ a complex number lying in the upper half plane. Note that by Lemma~\ref{lem:generallatprod} both $\Lambda_1,\Lambda_2$ have CM by an order in the imaginary quadratic number field $K=\Q(\tau_1)=\Q(\tau_2)$, and $\Lambda_1\Lambda_2$ is also a lattice with CM by an order in $K$.
    
    Let $\Gamma=\Lambda_1 \times \Lambda_2$ and $\eps_i=\frac{1}{\ol{\tau_i}-\tau_i}$. Letting $\{e_1,e_2,e_3,e_4\}$ be the basis of $\Gamma^*$ dual to the basis $\{\smat{1\\0},\smat{\tau_1\\0},\smat{0\\1},\smat{0\\\tau_2}\}$ of $\Gamma$, the map $H^1(X,\Z)\simeq \Gamma^* \simeq \Z^4 \to \C^2 \simeq \H^1(X,\O)$ is given by $$e_1 \mapsto \smat{-\epsilon_1\tau_1 \\ 0}, e_2 \mapsto \smat{\epsilon_1 \\ 0}, e_3 \mapsto \smat{0 \\ -\eps_2\tau_2}, e_4 \mapsto \smat{0 \\ \eps_2}$$  using the computation in (\ref{eqn:computation}).
    
    We may then calculate the image of the map $\bigwedge^2(\Z^4) \simeq \Z^6 \to \bigwedge^2(\C^2) \simeq \C$ based on the images of the basis elements $e_i \wedge e_j$ to be the additive subgroup $$\Lat(\Lambda_1 \times \Lambda_2)={\epsilon_1\epsilon_2}\lr{1,\tau_1,\tau_2,\tau_1\tau_2} =\epsilon_1\epsilon_2\Lambda_1\Lambda_2 \subseteq \C,$$ which is evidently homothetic to the lattice $\Lambda_1\Lambda_2$.
\end{proof}


\begin{example}\label{ex:ExE} 
If $E=\C/\lr{1,\tau}$ is an elliptic curve with complex multiplication, it follows that $\Con(E \times E)$ is isomorphic to the elliptic curve $\C/\Lambda'$ for $\Lambda'=\lr{1,\tau,\tau^2}=\lr{1,\tau}\cdot\lr{1,\tau}$. While $\Lambda'$ can be shown to be isogenous to $\Lambda$ (hence $\Con(E \times E)$ is isogenous to $E$), it is not usually the case that $\Lambda, \Lambda'$ are homothetic, meaning $\Con(E \times E)$ is not isomorphic to $E$.

As an example of this, consider $\Lambda=\lr{1,\frac{1+2i}{3}}$. Then $\Lambda'=\lr{1,\frac{1+2i}{3},\frac{-3+4i}{9}},$ which we can compute to be equal to $\lr{\frac{1}{3},\frac{2i}{9}}$, which by visual inspection is not homothetic to $\Lambda$. (We may also compute the $j$-invariants of $\Lambda$ and $\Lambda'$ to see that the corresponding elliptic curves are not isomorphic; this is done for the lattices $\lr{3,1+2i}$ and $\lr{3,2i}$, to which $\Lambda$ and $\Lambda'$ respectively are homothetic, in Example~\ref{ex:zsqrt-36} of the Appendix.)
\end{example}


\begin{lemma}\label{lem:prodnell} If $X=E_1 \times \dots \times E_n$ for isogenous CM elliptic curves $E_i=\C/\Lambda_i$, the lattice $\Lat(X)=\Lat(\Lambda_1 \times \dots \times \Lambda_n)$ is the Cartesian product $$\prod_{1 \le i < j \le n} \Lat(\Lambda_i \times \Lambda_j) \subseteq \C^{\binom{n}{2}}.$$
\end{lemma}
\begin{proof}
    As in the proof of Lemma~\ref{lem:prod2ell} we can assume $\Lambda_i=\lr{1,\tau_i}$ for $\tau_i$ and set $\Gamma=\Lambda_1 \times \dots \times \Lambda_n$ and $\eps_i=\frac{1}{\ol{\tau_i}-\tau_i}$. We also let $\{e_1,\dots,e_{2n}\}$ be the basis of $\Gamma^*$ dual to the specified basis for $\Gamma$. Using the computation (\ref{eqn:computation}) we see that the map $H^1(X,\Z) \simeq \Z^{2n} \to \H^1(X,\O)\simeq \C^n$ can be described by $$e_{2i-1} \mapsto -\eps_i\tau_iv_i,e_{2i} \mapsto \eps_iv_i$$ for $\{v_1,\dots,v_n\}$ the standard basis of $\C^n$. 

    The elements $\{v_i \wedge v_j : 1 \le i < j \le n\}$ form a basis of $\H^2(X,\O) \simeq \bigwedge^2\left(\C^n\right)\simeq \C^{\binom{n}{2}}$. We may then calculate the map $\H^2(X,\Z)\simeq \bigwedge^2(\Z^{2n})=\Z^{\binom{2n}{2}} \to \H^2(X,\O)$ on the generators $e_r \wedge e_s$ of $\H^2(X,\Z)$. Then using Lemma~\ref{lem:prod2ell} we can deduce that the image of $\H^2(X,\Z)$ is the Cartesian product  $$\prod_{1 \le i < j \le n} \eps_i\eps_j\Lambda_i \Lambda_j = \prod_{1 \le i < j \le n} \Lat(\Lambda_i \times \Lambda_j) \subseteq \C^{\binom{n}{2}}  $$
    
     As in Lemma~\ref{lem:prod2ell}, each lattice $\Lat(\Lambda_i \times \Lambda_j)=\eps_i\eps_j\Lambda_i\Lambda_j$ is homothetic to $\Lambda_i\Lambda_j$, which is a lattice with CM by an order in the field $K=\Q(\tau_i)$ by Lemma~\ref{lem:generallatprod}.
\end{proof}


\begin{lemma}\label{lem:definingBmX} In the setting of Lemma~\ref{lem:prodnell}, 
the lattice $\Lat_m(X)$ for $2 < m \le n$ is a Cartesian product  $$\prod_{1 \le i_1 < i_2 < \dots < i_m \le n}  \Gamma_{i_1, i_2, \dots, i_m},$$ where $\Gamma_{i_1, i_2, \dots, i_m}$ is a lattice in $\C$ homothetic to the lattice product $\Lambda_{i_1}\Lambda_{i_2} \dots \Lambda_{i_m}.$
\end{lemma}
\begin{proof} We define $\tau_i,\eps_i,e_i,v_i$ as in the proof of Lemma~\ref{lem:prodnell}. The elements $\{v_{i_1} \wedge v_{i_2} \wedge \dots \wedge v_{i_m}: 1 \le i_1 < i_2 < \dots < i_m \le n\}$ form a basis of $\H^m(X,\O) \simeq \bigwedge^m\left(\overline{T}\right)\simeq \C^{\binom{n}{m}}$\footnote{Recall $\ol{T}=\Hom_{\C-\rm{antilinear}}(\C,\C)$.}. We may then calculate the map $\H^m(X,\Z)\simeq \bigwedge^m(\Z^{2n})=\Z^{\binom{2n}{m}} \to \H^m(X,\O)$ on the generators $e_{j_1} \wedge e_{j_2} \wedge \dots \wedge e_{j_m}$ of $\H^m(X,\Z)$ to be the additive subgroup formed by the Cartesian product $$\prod_{1 \le i_1 < i_2 < \dots < i_m \le n} \eps_{i_1}\eps_{i_2} \dots \eps_{i_m} \Lambda_{i_1}\Lambda_{i_2} \dots \Lambda_{i_m} \subseteq \C^{\binom{n}{m}}.$$

We can apply Lemma~\ref{lem:generallatprod} repeatedly to see that each lattice product $\Gamma_{i_1,\dots,i_m}=\eps_{i_1}\eps_{i_2} \dots \eps_{i_m} \Lambda_{i_1}\Lambda_{i_2} \dots \Lambda_{i_m} $ is a lattice in $\C$ with CM by an order in $K=\Q(\tau_i)$.
\end{proof} 

\section{Structure of Brauer-Jacobians and complex abelian varieties}\label{sec:structureofJacobians}
In this section, we study the structure of Brauer-Jacobians for complex abelian varieties. As we show in Section \ref{sub:descBX}, a product of elliptic curves that are all isogenous to a fixed curve $E$ has a Brauer-Jacobian that can be calculated explicitly using the class group of the endomorphism algebra of $E$. Then, in Section \ref{sub:haveJacobians}, we show that the abelian varieties with Brauer-Jacobians are precisely those that have such a decomposition. The result is a complete description of $\Con(X)$ for any abelian variety of maximal Picard number.

\subsection{Computing the Brauer-Jacobian with class groups}\label{sub:descBX}

By Lemma~\ref{lem:prod2ell} we know that $\Lat(\Lambda_1 \times \Lambda_2)$ is homothetic to $\Lambda_1\Lambda_2$ for isogenous CM lattices $\Lambda_1,\Lambda_2$. We can say something further in the case where our pair of isogenous CM lattices have the same endomorphism ring. This ring $R$ must necessarily be an order in an imaginary quadratic number field; we write $\Cl(R)$ for its ideal class group.

\begin{definition} There is a bijective correspondence between the set of isomorphism classes of elliptic curves with CM by $R$ and $\Cl(R)$. Namely, if our elliptic curve $E$ is $\C/\Lambda$ for the lattice $\Lambda=\lr{1,\tau}$ and $\End(E)=\End(\Lambda)=R$, then $\Lambda \subseteq {\rm Frac}(R)$ must be an invertible fractional ideal of $R$. We will write $[E]$ or $[\Lambda]$ for the element of $\Cl(R)$ associated to $E=\C/\Lambda$ via this correspondence. \end{definition}

\begin{lemma}\label{lem:BE1E2sameendo} Let $R$ be an order in an imaginary quadratic number field $K$ with class group $\Cl(R)$. For $\Lambda_1, \Lambda_2$ invertible fractional ideals of $R$, the quotients $E_i = \C/\Lambda_i$ are CM elliptic curves with $\End(E_1)=\End(E_2)=R$.  Then $\Con(E_1 \times E_2)$ is a complex elliptic curve such that $$[E_1][E_2]=[\Con(E_1\times E_2)] \in \Cl(R).$$

\end{lemma}
\begin{proof}
By Lemma~\ref{lem:prod2ell} we know that $\Lat(\Lambda_1 \times \Lambda_2)$ is homothetic to $\Lambda_1 \Lambda_2$. By Lemma~\ref{lem:latprodidealprod} it is clear that $[\Lambda_1][\Lambda_2]=[\Lambda_1 \Lambda_2] =[\Lat(\Lambda_1 \times \Lambda_2)]\in \Cl(R)$, from which it follows that $[\Con(E_1 \times E_2)]=[E_1][E_2] \in \Cl(R)$.
\end{proof}

In the specific case of $E_1=E_2$, it follows from Lemma~\ref{lem:BE1E2sameendo} that $$[\Con(E \times E)]=[E]^2 \in \Cl(\End(E)).$$ We can use this to give an example of an elliptic curve $E$ such that the minimal field of definition of the elliptic curve $\Con(E \times E)$ is distinct from that of $E$.\footnote{It may be possible to define an abelian surface which is $\C$-isomorphic to $E \times E$ over a field which does not contain $j(E)$.}

\begin{example}\label{ex:qjezsqrt-36}
We recall from Example~\ref{ex:zsqrt-36zsqrt-9} that the $j$-invariant of the lattice $\lr{3,2i}$ (which has endomorphism ring $\Z[6i]$) is real and generates the degree $4$ extension $\Q(\sqrt[4]{12})$ of $\Q$; however, the $j$-invariant of $\lr{3,1+2i}$ (with the same endomorphism ring) is not real and generates the extension $\Q(i\sqrt[4]{12})/\Q$. Furthermore, $\Q(i\sqrt[4]{12}) \cap \R = \Q(\sqrt{3})$ and $j(\lr{3,2i}) \notin \Q(\sqrt{3})$, so $j(\lr{3,2i}) \notin \Q(j(\lr{3,1+2i}))$. 

If $E$ is the elliptic curve $\C/\lr{3,1+2i} \cong \C/\lr{1,\frac{1+2i}{3}}$ discussed in Example~\ref{ex:ExE}, we see that while $E$ is $\C$-isomorphic to an elliptic curve defined over $\Q(i\sqrt[4]{12})$, the same is not true of the elliptic curve $\Con(E \times E) \simeq_\C \C/\lr{3,2i}$, since $j(\Con(E \times E)) \notin \Q(i\sqrt[4]{12})$.
\end{example}

For two CM elliptic curves $E_1,E_2$ which are isogenous, we can deduce from the corresponding case for lattices, Lemma~\ref{lem:generallatprod}(i), that $\End(E_1)$ and $\End(E_2)$ are orders in the same imaginary quadratic field $K=\Q(\alpha)$. Assuming that $\alpha$ is chosen so $\OK=\Z[\alpha]$, then for some integers $a,b$ we have $\End(E_1)=\Z[a\alpha],\End(E_2)=\Z[b\alpha]$. We will use the map $\phi$ defined in Proposition~\ref{prop:defphi} in our general version of Lemma~\ref{lem:BE1E2sameendo}.

 For two isogenous CM elliptic curves $E_1,E_2$, there exists an imaginary quadratic number field $K$ and integers $a$ and $b$ such that $\OK=\Z[\alpha]$, $\End(E_1)=\Z[a\alpha]$, and $\End(E_2)=\Z[b\alpha]$. Let $c=\gcd(a,b)$.

\begin{theorem}\label{thm:BE1E2}
With the immediately preceding notation, $\Con(E_1 \times E_2)$ is a complex elliptic curve and we have $$\phi_{c,a,\alpha}([E_1])\phi_{c,b,\alpha}([E_2])=[\Con(E_1\times E_2)] \in \Cl(\Z[c\alpha]).$$
Thus, in the particular case where $\End(E_1)=\End(E_2)$ (so $a=b=c$), we have $[E_1][E_2]=[\Con(E_1 \times E_2)] \in \Cl(\Z[a\alpha]).$
\end{theorem}
\begin{proof}
As before, we assume $E_i=\C/\Lambda_i$ with $\Lambda_i$ an invertible fractional ideal of $\End(E_i)$. Via Lemma~\ref{lem:generallatprod}(ii) and the homothety of $\Lat(\Lambda_1 \times \Lambda_2)$ and $\Lambda_1  \Lambda_2$ we see that $\End(\Con(E_1 \times E_2))=\Z[c\alpha]$ (where $c={\rm gcd}(a,b)$ as in the statement of the theorem). 

Since $\Lat(\Lambda_1 \times \Lambda_2)$ is homothetic to $\Lambda_1  \Lambda_2$, it suffices to show that $\Lambda_1  \Lambda_2=( \Z[c\alpha]\cdot\Lambda_1 )  (  \Z[c\alpha]\cdot\Lambda_2)$ by the definition of $\phi$, which is computed in the proof of Lemma~\ref{lem:generallatprod}.
\end{proof}

Hence the case where $E_1,E_2$ have different endomorphism rings can be reduced to the case where $a \mid b$---i.e., where $\End(E_1) \supseteq \End(E_2)$---by considering the maps $\phi_{a,b,\alpha}$.

This allows us to describe $\Con(X)$ and $\Con_m(X)$ for $X$ a product of $n$ pairwise isogenous CM elliptic curves.

\begin{theorem}\label{thm:bvarnell} Let $X\cong E_1 \times \dots \times E_n$, where $E_i \cong \C/\Lambda_i$ are pairwise isogenous CM elliptic curves. Then $$\Con(X)=\prod_{1 \le i < j \le n} \Con(E_i \times E_j).$$ \end{theorem}
\begin{proof} We can deduce from Lemma~\ref{lem:prodnell} that $$\Lat(\Lambda_1 \times \dots \times \Lambda_n)=\prod_{1 \le i < j \le n} \Lat(\Lambda_i \times \Lambda_j) \subseteq \C^{\binom{n}{2}},$$ and the theorem follows. \end{proof}

\begin{theorem}\label{thm:computingBmX} For $X$ as in Theorem~\ref{thm:bvarnell} and $2 \le m \le n$, we have $$\Con_m(X) \cong \prod_{1 \le i_1 < i_2 < \dots < i_m \le n} \Con_m(E_{i_1} \times E_{i_2} \times\dots \times E_{i_m}).$$ In the case $m=n$, if $\End(E_i)=\Or_i=\Z[f_i\alpha]$ is the order of conductor $f_i$ in the imaginary quadratic number field $K=\Q \otimes \Or_i=\Q(\alpha)$ for each $i$, then we have $$\left[\Con_n(E_1 \times \dots \times E_n)\right]=\prod_{i=1}^n \phi_{d,f_i,\alpha}([E_i]) \in \Cl(\Or_0)$$ for $\Or_0=\Z[d\alpha]$ the order in $K$ of conductor $d=\gcd(f_1,\dots,f_n)$. 
\end{theorem}
\begin{proof}
    In the case $m=n$, we can use Lemma~\ref{lem:definingBmX} to compute that $\Con_n (E_1 \times \dots \times E_n) \cong \C/(\Lambda_1 \dots \Lambda_n)$ since the corresponding lattices are homothetic. By repeatedly applying Theorem~\ref{thm:BE1E2}, we can calculate the lattice product $$[\Lambda_1 \dots \Lambda_n] = \prod_{i=1}^n \phi_{d,f_i,\alpha}([E_i]) \in \Cl(\Or_0).$$

    For the general case, we again use Lemma~\ref{lem:definingBmX} to compute \begin{align*}\Con_m(X) &\cong \prod_{1 \le i_1 < i_2 < \dots < i_m \le n} \C/\left( \Lambda_{i_1}\Lambda_{i_2} \dots \Lambda_{i_m}\right)\\&\cong \prod_{1 \le i_1 < i_2 < \dots < i_m \le n} \Con_m(E_{i_1} \times E_{i_2} \times\dots \times E_{i_m}),\end{align*}
    via a homothety of lattices for each subset of indices $\{i_1,\dots,i_m\} \subseteq \{1, \dots, n\}$. 
\end{proof}

\begin{remark}\label{rmk:BmXextra}
    In the case $m=n$, we see that $\Con_m(E_1 \times \dots \times E_m) \cong \C/(\Lambda_1 \dots \Lambda_m)$, so in fact \begin{align*}\Con_m(E_1 \times \dots \times E_m) &\cong \Con(\Con_{m-1}(E_1 \times \dots \times E_{m-1}) \times E_m)\\
    &\cong\Con(\Con(\Con_{m-2}(E_1 \times \dots \times E_{m-2}) \times E_{m-1}) \times E_m)\\
    &\cong \dots\\
    &\cong\underbrace{\Con(\Con(\dots \Con(\Con(}_{m-1\text{ times}}E_1 \times E_2) \times E_3) \times \dots )\times E_m).\end{align*}
\end{remark}

\subsection{Which complex tori have Brauer-Jacobians?}\label{sub:haveJacobians}

We have been able to describe the Brauer-Jacobian of a product of pairwise isogenous CM elliptic curves in the previous sections. It has been shown (\cite{sm74} for the case of surfaces and \cite{schoen} for the general case) that a complex abelian variety which is isogenous to a product of isogenous CM elliptic curves must be isomorphic to such a product. Such abelian varieties are exactly those with maximal Picard number.

In the following theorem we use this fact to show, via our previous explicit computations, that the Brauer-Jacobians for all $2$-maximal abelian varieties are not just complex tori but abelian varieties. 

\main
\begin{proof}
It is immediately clear that (1) implies (2); the reverse implication is given by \cite[Satz 2.4]{schoen}. It is also immediate that (4) implies (3). 

Assuming (1), we may apply Theorem~\ref{thm:bvarnell} to see that $\Con(X)$ exists and is a product of complex elliptic curves, hence an abelian variety. Thus (1) implies (4).

Now assume (3). The fact that $\Con(X)=\H^2(X,\O)/\Lat(X)$ is a complex torus implies that $\Lat(X)$, the image of $\H^2(X,\Z)$ in $\H^2(X,\O) \simeq \C^{\binom{n}{2}}$, is a free abelian group of rank $n^2-n$. Since $\H^2(X,\Z)$ is a free abelian group of rank $\binom{2n}{2}=2n^2-n$, this means that the kernel of $\H^2(X,\Z) \to \H^2(X,\O)$ is isomorphic to $\Z^{n^2}$. 

By \cite[Chapter I.2, Theorem of Appell-Humbert]{mumabvar}, this kernel is isomorphic to the Neron-Severi group ${\rm NS}(X) \simeq \rm{Pic}(X)/\rm{Pic}^0(X)$, which is always a free abelian group of finite rank (see either \cite[Corollary 5.3.9]{blabvar} or \cite[Chapter IV.19, Corollary 2, Theorem 3]{mumabvar}). Furthermore, the rank of the Neron-Severi group is $n^2$ exactly in the situation of (2) \cite[Proposition 5.2.1, Corollary 5.3.8]{blabvar}. Therefore (3) implies (2). 

Then (2) implies (1) implies (4) implies (3) implies (2).  By \cite[Proposition 5.2.1, Corollary 5.3.8]{blabvar}, (5) is equivalent to (2). Thus all conditions are equivalent, as we desired to show.
\end{proof}

The proof of \cite[Satz 2.4]{schoen} is ``algebraic'', relying on the fact that for an order in an imaginary quadratic number field, any torsion-free module can be decomposed into a direct sum of rank $1$ modules \cite[Theorem 1.7]{bass}. However, this decomposition is not unique; it follows that for $X \simeq E_1 \times \dots \times E_n$, the $E_i$ are also not unique. This ambiguity of decomposition is discussed in \cite{sm74} for the case $n=2$ using binary quadratic forms, which we restate here.

\begin{theorem}\label{thm:surfacedecomp}
If $\A$ is a $2$-maximal complex abelian surface with $\A = \C^2/\Gamma$, then $\A \simeq \C/\OA \times \Con(\A)$ for $\OA=\{z \in \C : z\Gamma \subseteq \Gamma\}$. It then follows that \begin{enumerate}[(i)] \item $\A \simeq \A'$ for two $2$-maximal complex abelian surfaces if and only if $\OA = \Or(\A')$ and $\Con(\A) \simeq \Con(\A')$;
\item and for $E_i=\C/\Lambda_i,i=1,\dots,4$ pairwise isogenous CM elliptic curves with CM by the imaginary quadratic field $K$, we have $E_1 \times E_2 \simeq E_3 \times E_4$ if and only if $\End(E_1) \cap \End(E_2) = \End(E_3) \cap \End(E_4)$ and $\Lambda_1\Lambda_2$ is homothetic to $\Lambda_3\Lambda_4$.
\end{enumerate}
\end{theorem}
\begin{proof}
Condition (ii) is proved in \cite[Proposition 4.5]{sm74}.

If $\A=\C^2/\Gamma$ has a Brauer-Jacobian $\Con(\A)$, by Theorem~\ref{thm:main} we may assume without loss of generality that $\A=\C/\Lambda_1 \times \C/\Lambda_2$ for $\Lambda_1,\Lambda_2$ lattices in $\C$ with CM by orders $\Or_1,\Or_2$ in an imaginary quadratic number field $K$ (as $\Lambda_1,\Lambda_2$ are isogenous). It is then clear that $\OA=\Or_1 \cap \Or_2$, and $\OA$ is the order in $K$ with conductor $F=\lcm(f_1,f_2)$ for $f_i$ the conductors of $\Or_i$. Also set $f_0=\gcd(f_1,f_2)$ (and note that $F=f_1f_2/f_0$). 

Due to \cite[Bemerkung 2.6]{schoen} (and similarly to \cite[Lemma 8]{bf68}) we see that the lattices $\Lambda_1 \times \Lambda_2 \subseteq \C^2$ and  $\OA \times \Lambda_1\Lambda_2$ are isomorphic to one another as $\OA$-modules, which proves the main statement of the theorem, from which (i) follows. \end{proof}

\begin{remark}\label{rmk:primitivity} Using the notation of Theorem~\ref{thm:surfacedecomp}, we see that the \emph{degree of primitivity} $m=F/f_0$ 
is also an invariant of $\A$, and in fact if $\End(\Con(\A))=\Z[\alpha]$ we have $\OA=\Z[m\alpha]$. Thus a $2$-maximal complex abelian surface $\A$ is classified up to isomorphism by $\Con(\A)$ and $m$. We will refer to $\A$ such that $m=1$ as \emph{primitive}.
\end{remark}

The uniqueness of the decomposition in Theorem~\ref{thm:surfacedecomp} can also be derived from the general case \cite[Theorem 8]{bf68}, though in this case where $n=2$ we can derive it more directly as we have done here. The general case proved in \cite{bf68} will be useful when we consider the case $n=3$ in Section~\ref{sub:xsimeqbx}.

We can also use Theorem~\ref{thm:surfacedecomp} to apply the results of Section~\ref{sub:fieldofdef} concerning the minimal field of definition of elliptic curves to the following question. A 2-maximal primitive abelian surface $\mc{A}$ can always be decomposed as a product $\mc{A}=E_1 \times E_2$ of elliptic curves over the ring class field $L_\Or$ of $\Or=\End(E_1)=\End(E_2)=\End(\Con(\A))$; when is this possible over the smaller field $\Q(j(\Con(\A)))$?

\decompfieldofdef
\begin{proof} It follows from Theorem~\ref{thm:surfacedecomp}(ii) that the condition that $\End(\Con(\A))=\OA$ for $\A$ a $2$-maximal complex abelian surface is equivalent to the condition that for any $E_1,E_2$ with $\A \simeq E_1 \times E_2$ we have $\End(E_1)=\End(E_2)=\OA$. We set $\Or=\OA$ for ease of writing in this proof. 
In fact, for any $E_1,E_2$ with CM by $\Or$, we have $\A \simeq E_1 \times E_2$ if and only if $[E_1][E_2]=[\Con(\mc{A})] \in \Cl(\Or)$ by Theorems~\ref{thm:BE1E2} and ~\ref{thm:surfacedecomp}. We can then apply Lemma~\ref{lem:qje} to see that $\Q(j(E_1))=\Q(j(E_2))=\Q(j(\Con(\A)))$ if and only if $[E_1][E_2]^{-1},[E_1],$ and $[E_2]$ all have order $\le 2$. Since $\Cl(\Or)$ is abelian, this implies that $[\Con(\A)]=[E_1][E_2]$ has order $\le 2$ also. This completes one direction of the proof.

If $[\Con(\A)]$ has order $\le 2$, then $\Q(j(\Con(\A)))=\Q(j(\C/\Or))$, so we can set $E_1=\Con(\A),E_2=\C/\Or$ to get the desired decomposition.\end{proof}

\begin{remark}
Note that in Proposition~\ref{prop:decompfieldofdef}, the group element $[\Con(\A)] \in \Cl(\Or)$ has order $\le 2$ if and only if $j(\Con(\A))$ is real. (See \cite[Theorem 10.9, Exercise 11.1]{cox11}.)
\end{remark}

\begin{example}\label{ex:qjezsqrt-36pt2}
Recall from Example~\ref{ex:qjezsqrt-36} that the order $\Zsi=\Z[\sqrt{-36}]$ in the field $K=\Q(i)$ has ring class field $L_\Zsi=K(\sqrt[4]{12})$ and class group $\Cl(\Zsi)=\Z/4\Z$. The element $[\lr{3,1+2i}]$ has order 4 and $[\lr{3,1+2i}]^2=[\lr{3,2i}]$.

Let $\A$ be the abelian surface $\A = (\C/\lr{3,1+2i})^2$. Then $\Con(\A)\simeq \C/\lr{3,2i}$; note that $j(\lr{3,2i}) \notin \Q(j(\lr{3,1+2i}))=\Q(i\sqrt[4]{12})$ (see Example~\ref{ex:zsqrt-36}). It hence initially appears that $\A$ does not decompose as a product of elliptic curves over $\Q(j(\Con(\A)))$. 

However, by Theorem~\ref{thm:surfacedecomp}(ii), we see that $\A \cong \C/\Zsi \times \C/\lr{3,2i}$, and $\Q(j(\Zsi))=\Q(j(\lr{3,2i}))=\Q(\sqrt[4]{12})$. Therefore $\A$ can be decomposed as a product of elliptic curves over the minimal field of definition of $\Con(\A)$, as we would expect from Proposition~\ref{prop:decompfieldofdef} since the order of $[\Con(\A)]=[\lr{3,2i}]$ is $2$.
\end{example}

\begin{example}\label{ex:qjezsqrt-27pt2}
Recall from Example~\ref{ex:qjezsqrt-27} that the order $\Ztws=\Z[\sqrt{-27}]$ in the field $\Q(\zeta_3)$ has ring class field $L_{\Ztws}=\Q(\zeta_3,\sqrt[3]{2})$ and class group $\Cl(\Ztws)=\Z/3\Z$ with elements $1=[\Ztws],[\la_1]=[\lr{3,1+\s}],$ and$[\la_2]=[\lr{3,2+\s}]=[\lr{3,1+\s}]^2$. 

Let $\A'$ be the abelian surface $\A' = (\C/\la_1)^2$. Then $\Con(\A')\simeq \C/\lr{3,2+\s}$; similarly to the previous example, note that $j(\la_2) \notin \Q(j(\la_1))=\Q(\zeta_3\sqrt[3]{2})$ (see Example~\ref{ex:zsqrt-27}). 

By Theorem~\ref{thm:surfacedecomp}(ii), we see that the only other way to decompose $\A'$ as a product of two elliptic curves is as $\A' = (\C/\la_1)^2 \cong \C/\Ztws \times \C/\la_2,$ but $j(\Ztws) \not\in\Q(j(\la_2))=\Q(\zeta_3^2\sqrt[3]{2})$---as we would expect from Proposition~\ref{prop:decompfieldofdef}, since the order of $[\Con(\A')]=[\la_2]$ is $3$, we can't decompose $\A'$ as a product over $\Q(j(\Con(\A')))$. 
\end{example}

\subsection{When is an abelian variety isomorphic to its $m$-Jacobian?}\label{sub:xsimeqbx}

If $X$ is a $2$-maximal complex abelian variety, then by Theorem~\ref{thm:main} we can assume $X=E_1 \times \dots \times E_n$ is a product of $n$ elliptic curves which are CM and pairwise isogenous; by Theorem~\ref{thm:computingBmX} we see that $X$ is $m$-maximal for all $m \le \dim X$ and $\Con_m(X)$ is a product of $\binom{n}{m}$ pairwise isogenous CM elliptic curves. It is evident that the only $m$ between $2$ and $n$ for which it is possible that $X \simeq \Con_m(X)$ is $m=n-1$.

Before considering for which $X=E_1 \times E_2 \times \dots \times E_n$ we have $X \simeq \Con_{n-1}(X)$, we will wish to check which invariants of $X$ classify it up to isomorphism. Similar to Theorem~\ref{thm:surfacedecomp}, where we showed that a $2$-maximal complex abelian surface $\A$ is uniquely classified up to isomorphism by $\Con(A)$ and the degree of primitivity $m$ (Remark~\ref{rmk:primitivity}) these invariants will be $\Con_n(X)$ and $n-1$ integers. 

\begin{lemma}\label{lem:ndecomp}
Let $X=E_1 \times E_2 \times \dots \times E_n$ be a product of CM elliptic curves with each endomorphism ring $\End(E_i)$ equal to an order $\Or_i$ in a fixed imaginary quadratic number field $K$ (hence the elliptic curves are pairwise isogenous). If $f_i$ is the conductor of $\Or_i$ in $\OK=\Z[\alpha]$, let $d=\gcd(f_1,\dots,f_n)$ and $N=\lcm(f_1,\dots,f_n)$. Then there exist unique positive integers $r_2,\dots,r_{n-1}$ with $d=r_1 \mid r_2 \mid r_3 \mid \dots \mid r_{n-1} \mid r_n=N$ such that if $R_i$ is the order in $K$ of conductor $r_i$ we have an isomorphism $$X \cong \C/R_n \times \C/R_{n-1} \times \dots \times \C/R_2 \times \Con_n(X).$$ 
\end{lemma}
\begin{proof} We have already determined $R_1,R_n$ since $r_1=d,r_n=N$. It remains to determine the $r_i$ for $2 \le i \le n-1$.

We note that $[\Con_n(X)] = \prod_{i=1}^n \phi_{d,f_i,\alpha}([E_i]) \in \Cl(R_1)$ by Theorem~\ref{thm:computingBmX}. Thus by repeatedly applying Theorem~\ref{thm:surfacedecomp} to $X$, we can show that $$X \cong \C/R_n \times \C/R_{n-1} \times \dots \times \C/R_2 \times \Con_n(X)$$ for orders $R_2, \dots, R_{n-1}$ such that $R_n \subseteq R_{n-1} \subseteq \dots \subseteq R_2 \subseteq R_1$, and we can take $r_i$ to be the conductor of $R_i$. The uniqueness of this decomposition is a consequence of \cite[Theorem 8]{bf68}. Namely, given the endomorphism algebra $K$, the lattice $\Gamma$ and hence the complex torus $X=\C^n/\Gamma$ are determined up to isomorphism by the sequence of conductors $d=r_1,r_2,\dots,r_n=N$ and $\Con_n(X)$ (up to isomorphism/homothety of $\Con_n(X)$). The lemma follows.
\end{proof}
\begin{remark}\label{rmk:primitiven} 
    In the case $n=2$, we see that $r_1=d, r_2=N=md$. 
\end{remark}

We can now describe when $X \simeq \Con(X)$, or the case $n=3$. 

\begin{proposition}\label{prop:xsimbx}
For $X$ a $2$-maximal complex abelian threefold, we have $X \simeq \Con(X)$ if and only if $X \cong (\C/\Or)^3$ for some order $\Or$ in an imaginary quadratic number field $K$.
\end{proposition}
\begin{proof} Since $X$ is a $2$-maximal complex abelian threefold, by Theorem~\ref{thm:main} we have $X\cong E_1 \times E_2 \times E_3$ a product of three elliptic curves with $\End(E_i)=\Or_i$ for orders $\Or_i$ in an imaginary quadratic number field $K$ (hence the elliptic curves have CM and are pairwise isogenous).

As in Lemma \ref{lem:ndecomp}, let $f_i$ be the conductor of $\Or_i$ in $K$, and let $d=\gcd(f_1,f_2,f_3), N=\lcm(f_1,f_2,f_3)$. If $r=f_1f_2f_3/(dN)$ we can see easily that $d \mid r \mid N$. 

Applying Theorem~\ref{thm:surfacedecomp} twice, we find $$X=E_1 \times E_2 \times E_3 \simeq \C/R_3 \times \C/R_2 \times \Con_3(X),$$ where $R_1,R_2,R_3$ are the orders in $\OK$ with conductors $d,r=f_1f_2f_3/(dN),N$ respectively, and $\Con_3(X)$ is an elliptic curve with CM by $R_1$ such that $[\Con_3(X)]=\prod_{i=1}^3 \phi_{d,f_i,\alpha}([E_i]) \in \Cl(R_1)$. 

We then use Theorem~\ref{thm:bvarnell} to compute $\Con(X)$. To get a standard representative of the isomorphism class of $\Con(X)$ as in Lemma~\ref{lem:ndecomp}, we use Theorem~\ref{thm:surfacedecomp}. 
Thus we obtain $$\Con(X) \simeq \C/R_2 \times \Con_3(X) \times \Con_3(X) \simeq \C/R_2 \times \C/R_1 \times E,$$ for $E$ an elliptic curve with CM by $R_1$ such that $[E]=[\Con_3(X)]^2 \in \Cl(R_1)$.

Then by Lemma~\ref{lem:ndecomp}, we see that $X \simeq \Con(X)$ if and only if $E \simeq \Con_3(X), N=r,$ and $r=d$. The latter two equalities imply that $f_1=f_2=f_3$. This also means that all the orders $R_i,\Or_i$ are equal.

It then follows that $X \simeq \Con(X)$ if and only if $f_1=f_2=f_3$ and $[E]=[\Con_3(X)] \in \Cl(R_1)=\Cl(\Or_1)$; since then $[\Con_3(X)]=[E_1][E_2][E_3]$ and $[E]=[\Con_3(X)]^2$, it follows that $[\Con_3(X)]$ is the identity in $\Cl(\Or_1)$, so $\Con_3(X) \cong \C/\Or_1$. Because in this case, $X \cong \C/\Or_1 \times \C/\Or_1 \times \Con_3(X)$, we see that $X \simeq \Con(X)$ is equivalent to $X \cong (\C/\Or_1)^3$, where $\Or_1$ can be any order in an imaginary quadratic number field $K$.
\end{proof}

\begin{remark}\label{rmk:wholelattice} By Proposition~\ref{prop:xsimbx} and Lemma~\ref{lem:ndecomp}, we see that a complex abelian threefold $X$ is isomorphic to $\Con(X)$ if and only if $X \cong E^3$ for a complex elliptic curve $E$ which is isomorphic to $\C/\Or$ for an order $\Or$ in an imaginary quadratic number field $K$. 

It therefore makes sense to ask which elliptic curves $E$ satisfy this condition, or, in other words, which lattices $\Lambda$ in $\C$ are homothetic to an order $\Or$. In order for a lattice $\Lambda=\lr{1,\tau}$ to be an order, it must be a ring, for which it is sufficient to have $\tau^2 \in \Lambda$---in other words, $\tau^2+m\tau+n = 0 $ for some integers $m,n$, so $\tau$ is an algebraic integer of degree $2$. In fact, this is a necessary and sufficient condition for $\lr{1,\tau}$ to be an order in the imaginary quadratic number field $K=\Q(\tau)$.
\end{remark}

\begin{definition}\label{def:con2axn}
    The $\Con_2=\Con$ operator takes $2$-maximal complex abelian threefolds to $2$-maximal complex abelian threefolds; thus we can say that $\Con$ \emph{acts} on the set of isomorphism classes of $2$-maximal complex abelian threefolds.
\end{definition}

\begin{remark}
    In Proposition~\ref{prop:xsimbx} and Remark~\ref{rmk:wholelattice} we characterize the fixed points of the $\Con$ action. 
\end{remark}

\begin{proposition}\label{prop:conorbitfinite}
    The orbits of the $\Con$ action on $2$-maximal complex abelian threefolds are finite.
\end{proposition}
\begin{proof}
    Consider some $2$-maximal complex abelian threefold $X \cong E_1 \times E_2 \times E_3$, with $K,d,N,f_i$ defined as in Lemma~\ref{lem:ndecomp}. Furthermore, set $r,R_i$ as in the proof of Proposition~\ref{prop:xsimbx}. 

    It's computed in the proof of Proposition~\ref{prop:xsimbx} that $X \cong \C/R_3 \times \C/R_2 \times \Con_3(X)$ and $\Con(X) \simeq \C/R_2 \times \C/R_1 \times E$ for $E$ an elliptic curve with CM by $R_1$ such that $[E]=[\Con_3(X)]^2 \in \Cl(R_1)$. A very similar computation using Theorem~\ref{thm:bvarnell} and Theorem~\ref{thm:surfacedecomp} shows that $$\Con(\Con(X)) \cong \C/R_1 \times E \times E \cong (\C/R_1)^2 \times E'$$ for $E'$ an elliptic curve with CM by $R_1$ such that $[E']=[E]^2=[\Con_3(X)]^4 \in \Cl(R_1)$.

    In fact, for $k \ge 2$ we see that applying $\Con$ to $X$ a total of $k$ times is isomorphic to $(\C/R_1)^2 \times E_{(k)}$ where $E_{(k)}$ is an elliptic curve with CM by $R_1$ such that $$[E_{(k)}]=[E]^{\left(2^{k-1}\right)}=[\Con_3(X)]^{\left(2^k\right)} \in \Cl(R_1).$$ Since $\Cl(R_1)$ is finite, up to isomorphism there are only finitely many isomorphism classes $[E_{(k)}]$. Thus there are finitely many isomorphism classes of $2$-maximal abelian threefolds in the orbit of $X$ for the $\Con$ action.
\end{proof}

\begin{proposition}\label{prop:xsimconmx} For a $2$-maximal complex abelian variety $X$ of dimension $n \ge 3$, the ($n-1$)-Jacobian $\Con_{n-1}(X)$ is isomorphic to $X$ if and only if there exists a CM elliptic curve $E$ with $[E]^{n-2}=1 \in \Cl(\End(E))$ such that $X \cong (\C/\End(E))^{n-1} \times E$.
\end{proposition}
\begin{proof} Since $X$ is $2$-maximal we have $X=E_1 \times E_2 \times \dots \times E_n$ a product of CM elliptic curves with $\End(E_i)=\Or_i$ for orders $\Or_i$ in an imaginary quadratic number field $K$. 

Let $f_i$ be the conductor of $\Or_i$ and let $d=\gcd(f_1,\dots,f_n), N=\lcm(f_1,\dots,f_n)$. Then by 
Lemma~\ref{lem:ndecomp} there exist unique integers $r_1,\dots,r_n$ with $r_1=d \mid r_2 \mid \dots \mid r_{n-1} \mid r_n=N$ and an isomorphism $$X=E_1 \times \dots \times E_n \simeq \C/R_n \times \dots \times \C/R_2 \times \Con_n(X)$$ where $R_i$ is the order in $K$ with conductor $r_i$. 

Using Theorem~\ref{thm:computingBmX} and applying Theorem~\ref{thm:surfacedecomp} repeatedly, we see that $$\Con_{n-1}(X) \simeq \C/R_2 \times \left(\Con_n(X)\right)^{n-1}\simeq \C/R_2 \times (\C/R_1)^{n-2} \times E,$$ for $E$ an elliptic curve with CM by $R_1$ such that $[E]=[\Con_n(X)]^{n-1} \in \Cl(R_1)$.

By Lemma~\ref{lem:ndecomp}, we see that $X \simeq \Con_{n-1}(X)$ if and only if $E \simeq \Con_n(X), N=r_2,$ and $d=r_1=r_2=\dots=r_{n-1}$. The latter two equalities imply that $N=d$, so $f_1=f_2=\dots=f_n$. This also means that all the orders $R_i,\Or_i$ are all equal to each other.

It then follows that $X \simeq \Con_{n-1}(X)$ if and only if $f_1=\dots=f_n$ and $[E]=[\Con_n(X)] \in \Cl(R_1)=\Cl(\Or_1)$; since we also have $[E]=[\Con_n(X)]^{n-1}$, it follows that $[\Con_n(X)]^{n-2}=[E]^{n-2}=1 \in \Cl(\Or_1)$. Because in this case, $X \cong \left(\C/\Or_1\right)^{n-1} \times \Con_n(X)$, we see that $X \simeq \Con_{n-1}(X)$ is equivalent to $X \cong (\C/\Or_1)^{n-1} \times E$ for an elliptic curve $E$ with CM by $\Or_1$ and $[E]^{n-2}=1 \in \Cl(\Or_1)$, as we desired to show.
\end{proof}

\begin{remark}\label{rmk:conmfixed} By Proposition~\ref{prop:xsimconmx} and Lemma~\ref{lem:ndecomp} we see that $X=E_1 \times \dots \times E_n \cong \Con_{n-1}(X)$ if and only if $\End(\Con_n(X))=\End(E_1)=\dots=\End(E_n)$ and $[\Con_n(X)]^{n-2} \in \Cl(\End(\Con_n(X)))$ is the identity.
\end{remark}

\begin{definition}\label{def:conmaxn}
    Similarly to the $\Con$ action, we note that the $\Con_{n-1}$ operator takes $2$-maximal complex abelian $n$-folds to $2$-maximal complex abelian $n$-folds; thus we can say that $\Con_{n-1}$ \emph{acts} on the set of isomorphism classes of $2$-maximal complex abelian $n$-folds.
\end{definition}

\begin{proposition}\label{prop:conmorbitfinite}
    The orbits of the $\Con_{n-1}$ action on $2$-maximal complex abelian $n$-folds are finite.
\end{proposition}
\begin{proof}
    Consider some $2$-maximal complex abelian $n$-fold $X \cong E_1 \times E_2 \times \dots \times E_n$, with $\End(E_i)=\Or_i$ for orders $\Or_i$ in an imaginary quadratic number field $K$. 

If $f_i$ is the conductor of $\Or_i$, then by 
Lemma~\ref{lem:ndecomp} there exist unique integers $r_1,\dots,r_n$ with $$r_1=d=\gcd(f_1,\dots,f_n) \mid r_2 \mid \dots \mid r_{n-1} \mid r_n=N=\lcm(f_1,\dots,f_n)$$ and an isomorphism $$X=E_1 \times \dots \times E_n \simeq \C/R_n \times \dots \times \C/R_2 \times \Con_n(X)$$ where $R_i$ is the order in $K$ with conductor $r_i$. 

Using Theorem~\ref{thm:computingBmX} and applying Theorem~\ref{thm:surfacedecomp} repeatedly, we see that $$\Con_{n-1}(X) \simeq \C/R_2 \times \left(\Con_n(X)\right)^{n-1}\simeq \C/R_2 \times (\C/R_1)^{n-2} \times E,$$ for $E$ an elliptic curve with CM by $R_1$ such that $[E]=[\Con_n(X)]^{n-1} \in \Cl(R_1)$.
    
    Using Theorem~\ref{thm:computingBmX} and applying Theorem~\ref{thm:surfacedecomp} again, we obtain that $$\Con_{n-1}(\Con_{n-1}(X)) \cong \C/R_1 \times E^{n-1} \cong (\C/R_1)^{n-1} \times E_{(2)}$$ for $E_{(2)}$ an elliptic curve with CM by $R_1$ such that $[E_{(2)}]=[E]^{n-1}=[\Con_{n}(X)]^{\left((n-1)^2\right)} \in \Cl(R_1)$.

    We proceed by induction. Assume that the abelian $n$-fold we obtain by applying $\Con_{n-1}$ to $X$ a total of $k$ times is isomorphic to $(\C/R_1)^{n-1} \times E_{(k)}$, where $E_{(k)}$ is an elliptic curve with CM by $R_1$ such that $$[E_{(k)}]=[E]^{\left((n-1)^{k-1}\right)}=[\Con_n(X)]^{\left((n-1)^k\right)} \in \Cl(R_1).$$ (We showed this above for $k=2$.) Then by Theorem~\ref{thm:computingBmX} and Lemma~\ref{lem:ndecomp}, $$\Con_{n-1}((\C/R_1)^{n-1} \times E_{(k)}) \cong \C/R_1 \times E_{(k)}^{n-1} \cong(\C/R_1)^{n-1} \times E_{(k+1)},$$ where $E_{(k+1)}$ is an elliptic curve with CM by $R_1$ such that $[E_{(k+1)}]=[E_{(k)}]^{n-1}$. Using the inductive hypothesis, we then see that $$[E_{(k+1)}]=[E_{(k)}]^{n-1}=[\Con_n(X)]^{\left((n-1)^{k+1}\right)} \in \Cl(R_1).$$

   Since $\Cl(R_1)$ is finite, up to isomorphism there are only finitely many isomorphism classes $[E_{(k)}]$. Thus there are finitely many isomorphism classes of $2$-maximal abelian $n$-folds in the orbit of $X$ for the $\Con_{n-1}$ action.
\end{proof}

\section{Kummer varieties and singular K3 surfaces}

In this section, we discuss higher-weight Jacobians for Kummer varieties and singular K3 surfaces, which are related to Kummer varieties of $2$-maximal abelian varieties.

\begin{definition}\label{def:kummer}
    The \emph{Kummer variety\/} of an abelian variety $\mathcal{A}$, denoted $\Kum(\mathcal{A})$, is 
    quotient of the blowup at the $2$-torsion points $\widetilde{\mathcal{A}} \to \mathcal{A}$ by a lift $\tilde{\iota}$ of the involution (cf. \cite[Example 3.3]{huybrechts}). If $\mathcal{A}$ is an abelian surface, the Kummer surface $\mathcal{K}(\mathcal{A})$ is a K3 surface.
\end{definition}

\subsection{Singular K3 surfaces}\label{sec:K3}

In this section we consider specifically Brauer-Jacobians for complex abelian surfaces, Kummer surfaces, and K3 surfaces.

The Brauer-Jacobian of a 2-maximal abelian surface or a singular K3 surface $X$ is closely connected to the associated transcendental lattice $T_X$ and the quadratic intersection form $Q_X$ on this lattice. 

\begin{definition} For $X$ a 2-maximal abelian surface or a singular K3 surface, let $T_X$ be the \emph{transcendental lattice}, or the sublattice of transcendental cycles in $\H^2(X,\Z)$, and let $Q_X$ be the corresponding positive definite even integral binary quadratic form. 

A quadratic form $Q$ has an associated matrix $\left(\begin{smallmatrix} 2a & b \\ b& 2c \end{smallmatrix} \right)$ for integers $a,b,c$. The \emph{degree of primitivity} of $Q$ is $m=\gcd(a,b,c)$. If $m=1$ we say $Q$ is \emph{primitive}. We write $\tilde{Q}$ for $\frac{1}{m}Q$. 
\end{definition}

\begin{remark}\label{rmk:classisom}
    We will combine the classification of 2-maximal abelian surfaces \cite[Theorem 3.1]{sm74} or singular K3 surfaces \cite[Theorem 4]{si77} by positive definite even integral quadratic forms $Q_X$ and the isomorphism of the form class group and the ideal class group described in \cite[Theorem 7.7]{cox11}. Namely, if $\Or$ is the unique order in an imaginary quadratic field of discriminant $D$, the image of the class of a \emph{primitive} $Q(x,y)=ax^2+bxy+cy^2$ under the isomorphism from the form class group $\Cl(D)$ to $\Cl(\Or)$ is the class of the fractional ideal $\langle a, (-b+\sqrt{D})/2 \rangle$ of $\Or$. 
\end{remark}

\begin{proposition}\label{prop:quadformbrauerjac}
    Let $\A$ be a $2$-maximal abelian surface with ${\rm disc}(Q_\A)=D$ and ${\rm disc}(\widetilde{Q_\A})=D_0$, with $\Or_0$ the unique order in an imaginary quadratic field of discriminant $D_0$. 
    
    Then under the isomorphism $\Cl(D_0) \cong \Cl(\Or)$ of Remark~\ref{rmk:classisom} we have $[\widetilde{Q_\A}] = [\Con(\A)]$. Furthermore the degree of primitivity of $Q_\A$ is equal to the degree of primitivity of $\A$ defined in Remark~\ref{rmk:primitivity}. 
\end{proposition}
\begin{proof}
    This is proven in \cite[Corollary 3.3]{sm74}. 
    
    We can also see it explicitly from the construction in \cite[Equation 3.5]{sm74} that $\A \cong \C/\langle  a, (-b+\sqrt{D})/2 \rangle \times \C/\langle  1, (b+\sqrt{D})/2 \rangle$ and compute explicitly that $(b+\sqrt{D})/2$ is an algebraic integer (a root of $x^2-bx+ac$). 
    
    Then we see that $\langle  a, (-b+\sqrt{D})/2 \rangle$ is homothetic to the lattice associated to $\widetilde{Q_\A}$ by Remark~\ref{rmk:classisom}, and we can compute with Theorem~\ref{thm:BE1E2} that $\Con(\A) \cong \C/\langle  a, (-b+\sqrt{D})/2 \rangle$, so we obtain $[\widetilde{Q_\A}] = [\Con(\A)]$. The desired equality of the two notions of degree of primitivity follows from \cite[Equation 4.12, Proposition 4.5]{sm74}.
\end{proof}

\renewcommand{\thmkthreepream}{Let $Y$ be a singular K3 surface with ${\rm disc}(Q_Y)=D$ and ${\rm disc}(\widetilde{Q_Y})=D_0$, and let $\Or_0$ be the unique order in an imaginary quadratic field of discriminant $D_0$. }

\thmkthreebrauerjac
    
\begin{proof}
    In \cite{si77}, Shioda and Inose define an elliptic curve $C_X$ associated to a singular K3 or abelian surface $X$ which coincides with our definition of the Brauer-Jacobian $\Con(X)$. Furthermore, they proved (via a correspondence with equivalence classes of binary quadratic forms) that any singular complex K3 surface $Y$ (of maximal Picard number 20) is a double cover of a Kummer surface $X=\mathcal{K}(\mathcal{A})$ associated to a $2$-maximal abelian surface $\mathcal{A}$, with $Q_Y=\frac{1}{2}Q_X=Q_\A$ \cite[Theorem 4]{si77}. Then $\widetilde{Q_Y}=\widetilde{Q_X}=\widetilde{Q_\A}$, and $\Con(Y) \cong \Con(X) \cong \Con(\A)$ \cite[p.133-134]{si77}. 
\end{proof}

We now discuss fields of definition of singular K3 surfaces, including singular Kummer surfaces. 

\begin{proposition}\label{prop:k3fieldofdef}
    Let $X$ be a singular K3 surface with  intersection form $Q$ and degree of primitivity $m$. Set $\widetilde{Q}=\frac{1}{m}Q$. 
    
    If $[\widetilde{Q}]^2=1$ in the form class group $\Cl({\rm disc}(\widetilde{Q}))$, then $X$ has a model over the field $\Q(j(\Or))$ where $\Or$ is the unique imaginary quadratic order with discriminant $D={\rm disc}(Q)=m^2{\rm disc}(\widetilde{Q})$. 
\end{proposition}

 \begin{proof} Let $D_0=D/m^2={\rm disc}(\widetilde{Q})$ and, similarly to $\Or$, let $\Or_0$ be the unique imaginary quadratic order with discriminant $D_0$. We then have $\Or=1+m\Or_0$. 
 
Using the Shioda-Inose construction of \cite{si77}, Schuett proves that a singular K3 surface with intersection form $Q$ 
can be defined over the field $\Q(j(\Or),j(E))$ \cite[Proposition 10]{k3fieldofdef},  where $E$ is an elliptic curve associated to the isomorphism class $[\widetilde{Q}]$ in the form class group $\Cl(D_0) \cong \Cl(\Or_0)$ \cite[Theorem 7.7]{cox11} (see also Remark~\ref{rmk:classisom}).

Identifying $[\widetilde{Q}] \in \Cl(D_0)$ and $[E] \in \Cl(\Or_0)$, we obtain $[E]^2 = 1 = [\Or_0]^2 \in \Cl(\Or_0)$. Hence $\Q(j(\Or_0))=\Q(j(E))$ by Lemma~\ref{lem:qje}. By Lemma~\ref{lem:qjephi} we have $\Q(j(\Or_0)) \subseteq \Q(j(\Or))$, so $X$ has a model over $\Q(j(\Or),j(E)) = \Q(j(\Or))$ as we desired to show. 
\end{proof}

\kumfieldofdef
\begin{proof} Set $E=\Con(\A)$ for ease of notation. Using the Shioda-Inose construction of \cite{si77}, Schuett proves that the Kummer surface
$\mathcal{K}(\A)$ of a $2$-maximal abelian surface $\A$
can be defined over the field $\Q(j(\Or(\A)),j(E))$ \cite[Proposition 10, Cor. 12]{k3fieldofdef}. 

If $[E]$ has order $\le 2$ in $\Cl(\Or)$, then by Lemma~\ref{lem:qje} we have $\Q(j(\Or)) = \Q(j(E))$. By Lemma~\ref{lem:qjephi} we have $\Q(j(\Or)) \subseteq \Q(j(\Or(\A)))$, so $\mathcal{K}(\A)$ has a model over $  \Q(j(\OA))=\Q(j(\OA),j(E))$ as we desired to show. 
\end{proof}

\begin{remark}
    Note that if $\A$ has intersection form $Q$, then $\Kum(\A)$ has intersection form $2Q$; thus if $\A$ is has degree of primitivity $m_\A$, then $m_{\Kum(\A)}=2m_\A$. Thus for K3 surfaces with \emph{even} degree of primitivity and hence can be realized as a Kummer surface, Proposition~\ref{prop:kumfieldofdef} gives a model over a lower degree field than Proposition~\ref{prop:k3fieldofdef}.
\end{remark}

\begin{remark}\label{rmk:k3fieldofdef}
The minimal field of definition of a singular K3 surface associated to an intersection form $Q$ (or elliptic curve $C_X$) which is trivial in the class group, is discussed in \cite[Lemma 20]{k3fieldofdef}; the cases considered have degree of primitivity $m=1,2,$ or $3$ and certain discriminant values $D$. These are precisely the values of $m, D$ for which $\#\Cl(D)$ is both odd and equal to $\#\Cl(D/m^2)$, in which case the field of definition $\Q(j(\Or))$ is minimal \cite[Corollary 16, Lemma 19]{k3fieldofdef}.
\end{remark}

We end this section with an example of a 2-maximal abelian surface $\A$ which has a different minimal field of definition from its Kummer surface $\Kum(\A)$. 

\begin{example}\label{ex:akafieldofdef}
We recall from Example~\ref{ex:qjezsqrt-27}, that the order $\Ztws=\Z[\sqrt{-27}]$ in the field $\Q(\zeta_3)$ has ring class field $L_{\Ztws}=\Q(\zeta_3,\sqrt[3]{2})$ and class group $\Cl(\Ztws)=\Z/3\Z$, with the identity being $1=[\Ztws]$, and generator $[\la_1]=[\lr{3,1+\s}]$; we can compute $[\lr{3,1+\s}]^2=[\lr{3,2+\s}]=[\la_2]$. 

As computed in Example~\ref{ex:qjezsqrt-27pt2}, the elliptic curves $\C/\la_1, \C/\la_2 \cong \Con(\C/\la_1 \times \C/\la_1)$ have distinct minimal fields of definition $\Q(j(\la_1))=\Q(\zeta_3\sqrt[3]{2})$ and $\Q(j(\la_2))=\Q(\zeta_3^2\sqrt[3]{2})$. 

If $X=\Kum((\C/\la_1)^2)$, then by \cite[Corollary 16, Lemma 19]{k3fieldofdef}, the degree of any field of definition of $X$ must be divisible by $6$. In fact by \cite[Proposition 10, Cor. 12]{k3fieldofdef} there is a model for $X$ over the field $\Q(j(\Ztws), j(\la_1))=\Q(\zeta_3,\sqrt{3}[2])$ which has degree $6$ over $\Q$; this field of definition is minimal. 

However, $(\C/\la_1)^2$ clearly has a model over the field $\Q(\zeta_3\sqrt[3]{2})=\Q(j(\la_1))$, which has degree $3$ over $\Q$. It is then evident that the minimal degree for a field of definition of $(\C/\la_1)^2$ is less than the minimal degree for a field of definition of $X$.\end{example}

\subsection{Kummer varieties and generalities on higher-weight Jacobians}

\label{sec:generalities}
In this section, we briefly discuss higher-weight Jacobians in general, then discuss applications for blowups and equivariant $m$-Jacobians. Fix a smooth, projective, complex variety $X$. We assume $m > 0$ throughout this section. 

\begin{lemma}\label{lem:vtovijsurjective}
    Let $V$ be a pure real Hodge structure of weight $m$ with Hodge decomposition $V_{\C}=\bigoplus V^{i,j}$. For any $i$ and $j$ such that $i \ne j, i+j=m$, the $\R$-linear map
    $V\to V^{i,j}$ (from a real vector space to a complex vector space) is surjective. 
\end{lemma}
\begin{proof}
    Write $W^{i,j}=V^{i,j}\oplus V^{j,i}$, and let $W_0^{i,j}$ be the set of elements of the form $\alpha+\overline\alpha$ for $\alpha\in V^{i,j}$. By basic descent theory, the map $W_0^{i,j}\otimes\C\to W^{i,j}$ is an isomorphism. Moreover, the description of the invariants shows that the natural map $W_0^{i,j}\to V^{i,j}$ induced by the projection is an isomorphism of real vector spaces. Since $W_0^{i,j}\subset V$ (as $V$ is itself identified with the invariants of $V_{\C}$ under conjugation), the result follows.
\end{proof}

\begin{corollary}
    Given a pure integral Hodge structure $H$ of weight $m$, the subgroup $\Lat(H)\subset H^{0,m}$ is a lattice if and only if one (and hence both) of the following equivalent conditions hold:
    \begin{enumerate}
        \item $\rk\Lat(H)=2\dim_{\C} H^{0,m}$;
        \item $\dim_{\Q}\im(H\otimes\Q\to H^{0,m})=2\dim_{\C} H^{0,m}$.
    \end{enumerate}
\end{corollary}
\begin{proof}
    By Lemma \ref{lem:vtovijsurjective}, the map $\Lat(H)\otimes\R\to H^{0,m}$ is a surjection, which shows that $\Lat(H)$ spans the real vector space underlying $H^{0,m}$.  Thus we can choose a saturated sublattice $\Lambda\subset\Lat(H)$ that is a real basis for $H^{0,m}$, yielding an identification $\Lambda\otimes\R\to H^{0,m}$. The resulting topological group (torus) $G:=H^{0,m}/\Lambda$ is compact, so any element $\gamma\in\Lat(H)\setminus\Lambda$ will have multiples in $G$ that cluster, showing that $\Lat(H)$ cannot be discrete in $H^{0,m}$.

    The equivalence of the two conditions follows from the fact that the rank of $\Lat(H)$ is equal to the dimension of its $\Q$-span inside of $H^{0,m}$.
\end{proof}
Thus there is a purely numerical criterion for $\Lat(H)$ to be discrete in $H^{0,m}$. This will be useful for studying the properties of higher-weight Jacobians under natural geometric operations.

\begin{definition}
    Given a pure integral Hodge structure $H$ of weight $m$, the \emph{Jacobian discrepancy\/} of $H$ is the number
    $\delta(H)=\rk\Lat(H)-2\dim_{\C} H^{0,m}$ (see Notation \ref{notation:latcon}). If $H=\H^m(X,\Z)$, we will write $\delta_m(X)$ for $\delta(H)$.
\end{definition}
\begin{observation}
    We can also read $\delta(H)$ off from the rational Hodge structure $H_{\Q}:=H\otimes_{\Z}\Q$: writing $L(H_{\Q})$ for the image of the map of rational vector spaces $H_{\Q}\to H^{0,m}$, we have $\delta(H)=\dim_{\Q}L(H_{\Q})-2\dim_{\C}H^{0,m}.$
\end{observation}

\begin{remark}\label{rmk:discrepancyisadditive} By Lemma~\ref{lem:vtovijsurjective}, the map $\Lat(H)\otimes\R\to H^{0,m}$ is a surjection, so $\rk \Lat(H) \ge 2\dim_\C H^{0,m}$; hence $\delta(H)$ is always nonnegative.
By definition, $H$ has a Jacobian if and only if $\delta(H)=0$. Given a collection of Hodge structures $H_i$, we see that $\delta(\bigoplus H_i)=\sum\delta(H_i)$. 
\end{remark}
\begin{lemma}\label{lem:Jacobiantorsiondimension} 
    Let $H$ be a pure Hodge structure of weight $m$ and let $r=\dim_{\C} H^{0,m}$. For any positive prime number $p$, we have 
    $$\dim_{\F_p}\Con(H)[p]=2r+\delta(H).$$
\end{lemma}
\begin{proof}
    Since $\Lat(H)$ spans $H^{0,m}$, we can choose a saturated lattice $\Lambda\subset\Lat(H)\subset H^{0,m}$. It follows that $\Con(H)$ is isomorphic to $(\R/\Z)^{2r}/\Gamma$ for some finitely generated free abelian subgroup $\Gamma$ of non-negative rank equal to $\delta(H)$. By dividing elements of $\Gamma$ by $p$, it follows that we have $\dim_{\F_p}\Con(H)[p]=2r+\delta(H)$, as desired.
\end{proof}

\begin{corollary}\label{cor:addjac}
    If $H_i$, $i=1,\ldots,n$, are pure Hodge structures of weight $m$ then $H=\bigoplus H_i$ has a Jacobian if and only if each $H_i$ has a Jacobian.
\end{corollary}
\begin{proof}
    By construction, we have $\Con(H)=\bigoplus \Con(H_i)$, so the result follows immediately from Remark~\ref{rmk:discrepancyisadditive}.
\end{proof}

\begin{corollary}\label{cor:projbundle}
    If $P\to X$ is a projective bundle over $X$ of relative dimension $d$, then $P$ has an $m$-Jacobian if and only if $X$ has an $m$-Jacobian.
\end{corollary}
\begin{proof}
    The projective bundle formula gives an isomorphism of Hodge structures
    $$\H^m(P,\Z)\cong\bigoplus_{i=0}^d\H^{m-2i}(X,\Z)(-i).$$
    If $H_i$ is the weight $m$ Hodge structure $\H^{m-2i}(X,\Z)(-i)$, we note that $H_i^{0,m}$ is zero for $i > 0$; thus $H_i$ has a trivial $m$-Jacobian for $i > 0$. The result then follows from Corollary \ref{cor:addjac}.
\end{proof}

\begin{proposition}\label{prop:blow up Jacobian}
    Suppose $Y\subset X$ is a smooth proper subvariety of codimension $d$ and $\pi:\widetilde X\to X$ is the blowup along $Y$. Then
    $\delta_m(\widetilde X)=\delta_m(X)$.
    Thus $\widetilde X$ has an $m$-Jacobian if and only if $X$ has an $m$-Jacobian.

    Moreover, if $\widetilde X$ has an $m$-Jacobian then the natural map $$\Con_m(X)\to\Con_m(\widetilde X)$$
    is an isomorphism.
\end{proposition}
\begin{proof}
   Let $E$ be the preimage of $Y$ in $\widetilde X$. By Theorem 7.31 of \cite{voisinthing}, there is an isomorphism of Hodge structures
   $$\H^m(\widetilde X,\Z)\cong\H^m(X,\Z)\oplus\bigoplus_{i=1}^{d-1}\H^{m-2i}(Y,\Z)(-i).$$ As in the proof of Corollary~\ref{cor:projbundle}, we note that the weight $m$ Hodge structures $H_i:=\H^{m-2i}(Y,\Z)(-i)$ have $H_i^{0,m}=0$. It follows that $\delta_m(\widetilde X)=\delta_m(X)+\sum_{i=1}^{d-1} \delta(H_i)=\delta_m(X)$.

    To see the final assertion, note that for a blowup along a smooth center, the higher derived functors satisfy $R^i\pi_*\O_{\widetilde X}=0$ (since smooth varieties have rational singularities). Hence the natural map $\H^m(X,\ms O_X)\to\H^m(\widetilde X,\ms O_{\widetilde X})$ is an isomorphism that induces a map of $m$-Jacobians $\Con_m(X)\to\Con_m(\widetilde X)$. The result follows.
\end{proof}

\begin{corollary}\label{cor:birationalinvariance}
    Let $X$ and $X'$ be birational smooth projective varieties. Then $X$ has a $q$-Jacobian if and only if $X'$ has a $q$-Jacobian, for any $q$.
\end{corollary}
\begin{proof}
    By the Weak Factorization Theorem, we can compare $X$ and $X'$ by a chain of blowups and blowdowns with smooth centers. Following the diagram of these morphisms and using Proposition \ref{prop:blow up Jacobian} gives the desired result.
\end{proof}

Now suppose $H$ is a pure polarizable Hodge structure of weight $m$ admitting an action of a finite group $G$ (in the category of Hodge structures).

\begin{lemma}\label{lem:invariant}
    If $H$ has a Jacobian then the invariant substructure $H^G$ has a Jacobian; moreover, $\Con(H)$ admits a $G$-action and the invariant connected subtorus $(\Con(H)^G)^0$ is isogenous to $\Con(H^G)$. 
\end{lemma}
\begin{proof}
    By flatness, we have $(H^G)_{\Q}=(H_{\Q})^G$. Thus it suffices to consider the rational Hodge structure $H_{\Q}$. By semisimplicity, we have a decomposition $H_{\Q}=H_{\Q}^G\oplus H'$ of Hodge structures. 

    Then by Remark~\ref{rmk:discrepancyisadditive}, $\delta(H^G)=0$, so $H^G$ has a Jacobian.

    The Jacobian $\Con(H)$ naturally inherits the $G$-action of $H$. However, because of the distinction between invariant tangent vectors and fixed points of a group acting on complex torus, when finding $\Con(H^G)$ we must consider specifically the identity component $(\Con(H)^G)^0$: the invariant connected subtorus of $\Con(H)$. Furthermore, unlike Corollary \ref{cor:addjac}, we have a decomposition of rational Hodge structures rather than integral Hodge structures, so we may only conclude that the natural map $\Con(H^G) \to (\Con(H)^G)^0$ is an isogeny.
\end{proof}

\begin{proposition}\label{prop:quotientjac}
    Let $X$ be a smooth projective variety with an action by a finite group $G$ such that the quotient $Y=X/G$ is smooth. If $X$ has a $q$-Jacobian then $Y$ has a $q$-Jacobian. Morever, $\Con_q(X)$ admits a $G$-action and $(\Con_q(X)^G)^0$ is isogenous to $\Con(Y)$. 
\end{proposition}
\begin{proof}
    By the Cartan-Leray spectral sequence (Theorem XVI.8.4 of \cite{cartaneilenberg}), we have $$\H^q(Y,\Q)=\H^q(X,\Q)^G.$$ The result thus follows from Lemma \ref{lem:invariant}.
\end{proof}

We are now able to (partially) generalize Theorem~\ref{thm:k3brauerjac}:
\kummercon
\begin{proof}
    By Proposition \ref{prop:blow up Jacobian}, the blow up $\widetilde{\mathcal{A}}$ of $\mathcal{A}$ at the $2$-torsion points has a $2m$-Jacobian with $\Con_{2m}(\A) \cong \Con_{2m}(\widetilde{\A})$. 
    
    The involution $\iota$ has a trivial action on the even degree cohomology of $\A$, so the invariant connected subtorus of $\Con_{2m}(\A) $ is the whole torus. Then by Proposition \ref{prop:quotientjac}, the Kummer variety $\Kum(\mathcal{A})$ has a $2m$-Jacobian which is isogenous to $\Con_{2m}(\widetilde{\A}) \cong \Con_{2m}(\A)$. Moreover, for any $q$ odd we have $\Con_{q}(\Kum(\A))=0$ by Proposition~\ref{prop:quotientjac}.
\end{proof}

\appendix
\section{Calculating $j$-invariants}\label{app}

To compute the class group of an order and the homothety classes of lattices with CM by that order, we used Sage to list all primitive positive definite forms of a given discriminant and compute the structure of the form class group.

Then by the bijection between the form class group and the ideal class group described in \cite[Theorem 7.7]{cox11} (see Remark~\ref{rmk:classisom}) we can then compute the homothety classes of lattices with CM by $\Or_D$.

To compute the ring class field $L_\Or$, we can adjoin a root of the Hilbert class polynomial $H_\Or(X)$ to $K={\rm Frac}(\Or)$ \cite[Proposition 13.2]{cox11}. We used Sage's database of Hilbert class polynomials and found a simpler description of $L_\Or$, then factored $H_\Or(X)$ over $L_\Or$ to get exact values of the the $j$-invariants of the lattices with CM by $\Or$. To identify which $j$-invariant corresponds to which lattice, we need only then approximate the $j$-invariants of each lattice; this can be done by truncating the $q$-series described in \cite[Proposition 9.12]{wa08}.

Thus we may compute the following, using $\sim$ to denote homothety of lattices. Below we fix a cube root of unity $\zeta_3=\frac{-1+\s}{2}$.

\begin{example}\label{ex:zsqrt-9}
The order $\Z[3i]=\Z[\sqrt{-9}]$ in $\Q(i)$ has discriminant $-36$; the class group is isomorphic to $\Z/2\Z$, with the two classes being represented by the quadratic forms $x^2+9y^2,2x^2+2xy+5y^2$ or the lattices $\lr{1,3i}$ (identity) and $\lr{2,-1+3i} \sim \lr{3,1+i}$, with $j$-invariants
\begin{align*}
    j(\lr{1,3i})&=76771008+44330496\sqrt{3},\\
    j(\lr{3,1+i})&=76771008-44330496\sqrt{3}.
\end{align*}
The ring class field is $L_{\Z[3i]}=\Q(i,\sqrt{3})$.
\end{example}
\begin{example}\label{ex:zsqrt-36}
The order $\Z[6i]=\Z[\sqrt{-36}]$ in $\Q(i)$ has discriminant $-144$; the class group is isomorphic to $\Z/4\Z$, with the four classes being represented by the quadratic forms (and corresponding lattices) given below with their orders: \begin{align*}x^2+36y^2: & \lr{1,6i},  &\text{(identity)}
\\5x^2\pm4xy+8y^2:& \lr{5,\pm2+6i} \sim \lr{3,1\pm2i},&\text{(order 4)}\\4x^2+9y^2: &\lr{4,6i}\sim\lr{3,2i}. &\text{(order 2)}\end{align*}  The $j$-invariant of the class $[\langle 3, 1+2i \rangle]^k$ for $k=0,1,2,3$ is 
\begin{align*}
j([\lr {3,1+2i}]) &=    5894625992142600+i^{2k}\cdot 3403263903336192\sqrt{3}\\& +i^k \cdot 3167093925247392\sqrt[4]{12}+i^{3k} \cdot   914261265145368\left(\sqrt[4]{12}\right)^3.
\end{align*}
 

The ring class field is $L_{\Z[6i]}=\Q(i,\sqrt[4]{12})$.
\end{example}

\begin{example}\label{ex:zsqrt-27}
The order $\Z[3\sqrt{-3}]=\Z[\sqrt{-27}]$ in $\Q(\sqrt{-3})$ has discriminant $-108$; the class group is isomorphic to $\Z/3\Z$, with the three classes being represented by the quadratic forms (and corresponding lattices) given below with their orders: \begin{align*}x^2+27y^2: & \lr{1,3\s},  &\text{(identity)} \\4x^2-2xy+7y^2&: \lr{4,1+3\sqrt{-3}} \sim \lr{3,1+\sqrt{-3}},  &\text{(order 3)}\\4x^2+2xy+7y^2&: \lr{4,-1+3\sqrt{-3}} \sim \lr{3,2+\sqrt{-3}}  &\text{(order 3)}\end{align*} The $j$-invariant of the class $[\langle 3,1+\sqrt{-3} \rangle]^k$ for $k=0,1,2$ is 
\begin{align*}
j([\lr{3,1+\sqrt{-3}}[^k)&=31710790944000\left(\zeta_3^k\sqrt[3]{2}\right)^2 + 39953093016000\zeta_3^k\sqrt[3]{2} + 50337742902000,
\end{align*}

The ring class field is $L_{\Z[3\s]}=\Q(\s,\sqrt[3]{2})$.
\end{example}

\begin{example}\label{ex:zsqrt-48}
The order $\Z[4\sqrt{-3}]=\Z[\sqrt{-48}]$ in $\Q(\sqrt{-3})$ has discriminant $-192$; the class group is isomorphic to $\Z/2\Z \times \Z/2\Z$, with the four classes being represented by the quadratic forms (and corresponding lattices) given below with their orders:\begin{align*}x^2+48y^2:&\lr{1,4\s},  &\text{(identity)}
\\7x^2+2xy+7y^2:&\lr{7,-1+4\s} \sim \lr{4,2+\s}, &\text{(order 2)}
\\4x^2+4xy+13y^2:&\lr{4,-2+4\s}\sim\lr{2,1+2\s},&\text{(order 2)}
\\3x^2+16y^2:&\lr{3,4\s}\sim\lr{4,\s}.&\text{(order 2)}
\end{align*} with $j$-invariants
\begin{align*}
j(\lr{1,4\s})&=820762881440077125\sqrt{6} + 1160733998424384000\sqrt{3} \\&
+ 1421603011620136125\sqrt{2} + 2010450259344609000,\\
j(\lr{4,2+\s})&=-820762881440077125\sqrt{6} - 1160733998424384000\sqrt{3} \\&
+ 1421603011620136125\sqrt{2} + 2010450259344609000, \\
j(\lr{2,1+2\s})&=-820762881440077125\sqrt{6} + 1160733998424384000\sqrt{3} \\&
- 1421603011620136125\sqrt{2} + 2010450259344609000, \\
j(\lr{4,\s})&=820762881440077125\sqrt{6} - 1160733998424384000\sqrt{3} \\&
- 1421603011620136125\sqrt{2} + 2010450259344609000. 
\end{align*}
The ring class field is $L_{\Z[4\s]}=\Q(i,\sqrt{2},\sqrt{3})$.
\end{example}





\subsection*{Declaration of generative AI}
During the preparation of this work the authors used ChatGPT Pro in order to verify the computations in the appendix and to review possible mathematical errors. After using this tool, the authors reviewed and edited the content as needed and take full responsibility for the content of the published article.

\printbibliography
\end{document}